\documentclass[11pt]{article}
\usepackage{amsthm, amsmath, amssymb, amsfonts, url, booktabs, tikz, setspace, fancyhdr, bm}
\usepackage{geometry}
\usepackage{enumerate}
\usepackage{algorithm}
\usepackage{algorithmic}
\usepackage[shortlabels]{enumitem}
\usepackage[babel]{microtype}
\usepackage[english]{babel}
\usepackage{comment}
\usepackage{bbm}
\usepackage{csquotes}
\usepackage{mathabx}
\usetikzlibrary{calc, angles, quotes, intersections}
\usepackage{graphicx}
\usepackage{float}
\usepackage{xcolor}
\counterwithin{figure}{section}
\newtheorem{theorem}{Theorem}[section]
\newtheorem{prop}[theorem]{Proposition}
\newtheorem{lemma}[theorem]{Lemma}
\newtheorem{cor}[theorem]{Corollary}

\theoremstyle{definition}

\newtheorem*{defn-non}{Definition}

\usepackage{hyperref}
\usepackage{cleveref}
\crefname{lemma}{lemma}{lemmas}
\Crefname{lemma}{Lemma}{Lemmas}
\crefname{prop}{proposition}{propositions}
\Crefname{prop}{Proposition}{Propositions}
\crefname{cor}{corollary}{corollaries}
\Crefname{cor}{Corollary}{Corollaries}
\crefname{claim}{claim}{claims}
\Crefname{claim}{Claim}{Claims}
\crefname{defn}{definition}{definitions}
\Crefname{defn}{Definition}{Definitions}
\newlist{Case}{enumerate}{2}
\setlist[Case, 1]{%
    label           =   {\bfseries Case \arabic*.},
    labelindent=1em, labelwidth=1.3cm, labelsep*=1em, leftmargin=!
}
\setlist[Case, 2]{%
    label           =   {\bfseries Subcase \arabic{Casei}.\arabic*.},
    labelindent=-1em, labelwidth=1.3cm, labelsep*=1em, leftmargin=!
}

\newcommand{\ceil}[1]{\left\lceil #1\right\rceil}
\newcommand{\floor}[1]{\left\lfloor #1\right\rfloor}

\newcommand{\ex}{\mathrm{ex}}
\newcommand{\cP}{\mathcal{P}}

\newcommand{\cI}{\mathcal{I}}
\newcommand{\cJ}{\mathcal{J}}

\newcommand*{\abs}[1]{\lvert#1\rvert}

\newcommand{\R}{\mathbb{R}}

\newcommand{\aff}{\operatorname{aff}}

\newcommand{\cF}{\mathcal{F}}

\newcommand{\cC}{\mathcal{C}}
\newcommand{\cL}{\mathcal{L}}

\newcommand{\E}{\mathbb{E}}

\newcommand{\Z}{\mathbb{Z}}
\newcommand{\cH}{\mathcal{H}}
\newcommand{\cT}{\mathcal{T}}
\newcommand{\HD}{\operatorname{HD}}

\newcommand{\rank}{\operatorname{rank}}

\numberwithin{equation}{section}

\title{Hadwiger--Debrunner numbers from rainbow Tur\'an numbers and beyond}
\author{Wenchong Chen\thanks{Department of Mathematics, Emory University, Atlanta, GA 30322, United States of America.
Email: wenchong.chen@emory.edu.}
\and
Gennian Ge\thanks{School of Mathematical Sciences, Capital Normal University, Beijing 100048, China. Email: gnge@zju.edu.cn. Gennian Ge is supported by the National Key Research and Development Program of China under Grant 2025YFC3409900, the National Natural Science Foundation of China under Grant 12231014, and Beijing Scholars Program.}
\and 
Yang Shu\thanks{School of Mathematical Sciences, University of Science and Technology of China, Hefei, 230026, China. Email: shuyangyyyy@mail.ustc.edu.cn.}
}
\date{}

\begin{document}
\maketitle
\begin{abstract}
The Hadwiger-Debrunner number $\HD_d(p,q)$ is defined as the maximum piercing number over all finite families of compact convex sets in $\mathbb R^d$ in which among every $p$ sets, some $q$ have non-empty intersection. 
We develop three complementary constructions that improve lower bounds for Hadwiger–Debrunner numbers. 
For any fixed $q>d\ge 2$, a randomly sparsified family of hyperplanes, analyzed using hypergraph containers, yields
$$\HD_d(p,q)\ge p^{(d+q-2)/(q-1)+o(1)}$$ as $p\rightarrow\infty$.
To deduce further bounds, we establish a transfer argument from upper bounds on the sizes of sufficiently general-positioned subsets of high-dimensional integer grids to lower bounds for Hadwiger–Debrunner numbers. 
Two methodologies are adopted to derive the former: 
one connects rainbow Tur\'an theory and general-positioned grid sets, 
and the other uses full spark linear reduction of tuples with bounded-size preimages. 
As a combinatorial ingredient, we describe the rainbow Tur\'an behavior under suspension operations, which may be of independent interest. 
\end{abstract}
\begin{section}{Introduction}

As a central result of combinatorial convexity, 
\emph{Helly's Theorem} \cite{Helly1923} asserts that, 
for a finite family $\cF$ of convex sets in $\mathbb R^d$, 
if any $d+1$ sets of $\cF$ have non-empty intersection, 
then all the sets in $\cF$ have non-empty intersection. 
As a natural yet challenging extension, 
Hadwiger and Debrunner relaxed the $(d+1)$-wise intersecting condition to the so-called \emph{$(p,q)$-property}, 
which requires that among every $p$ members of a family, some $q$ have a common point. 
The celebrated \emph{$(p,q)$-theorem} established by Alon and Kleitman \cite{AlonKleitman1992} states that, 
for $p\ge q\ge d+1$, every finite family of nonempty compact convex sets with $(p, q)$-property
can be pierced by a constant number of points, depending only on $d,p,q$. 
Numerous generalizations of these results have been explored in recent decades, 
encompassing colorful, fractional, quantitative and topological variants \cite{AmentaDeLoeraSoberon2017,Barany1982,BaranyKatchalskiPach1982,KatchalskiLiu1979}, 
as well as Helly-type and $(p,q)$-type results for special geometric families \cite{CarmiKatzMorin2023,DanzerGrunbaum1982,KellerSmorodinsky2020,KimNakprasitPelsmajerSkokan2006,NaszodiTaschuk2010,Tardos1995}.

For general compact convex families, 
one naturally asks what universal bound on the number of piercing points is guaranteed by $(p,q)$-theorem. 
This quantity is called the \emph{Hadwiger-Debrunner number}, 
denoted by $\HD_d(p,q)$, 
formally defined as the minimum integer $N=N(d,p,q)$, 
such that every family $\cF$ of compact convex sets in $\mathbb R^d$ with $(p,q)$-property 
admits a point set $P\subseteq \mathbb R^d$ satisfying $|P|\le N$ and $P\cap F\ne\emptyset$ for each $F\in\cF$. 
Substantial progress has been made on the bounds of $\HD_d(p,q)$. 
For upper bounds, Keller, Smorodinsky and Tardos \cite{KellerSmorodinskyTardos2018} concluded that, for any fixed $d,q$, $\mathrm{HD}_d(p,q)\le p^{d(q-1)/(q-d)+o(1)}$ when $p$ tends to infinity. 
For lower bounds, Keller and Smorodinsky \cite{KellerSmorodinskyLower} established the estimate $\mathrm{HD}_2(p,q)\ge p^{1+1/(4q-7)-o(1)}$ for every $q\ge 3$. 
More recently, Roche-Newton \cite{RocheNewton2026} improved the planar bound to $\mathrm{HD}_2(p,3)\ge p^{5/4-o(1)}$, 
while in higher dimensions, Keller and Smorodinsky \cite{KellerSmorodinskyUnions2026} proved that $\mathrm{HD}_d(p,d+1)\ge p^{2-o(1)}$ 
for $d\ge 3$ and a series of $p$ tending to infinity.

In this paper, we develop connections between Hadwiger-Debrunner numbers and various areas, and establish multiple constructions that improve the lower bounds. 
Our first construction works directly with randomly sparsified hyperplanes, where the absences of concurrences are controlled by hypergraph containers.
Additionally, we use point-hyperplane duality to reduce the lower bounds of HD numbers to the restrictions on the size of sufficiently general-positioned grid sets. 
To obtain the latter, we utilize two distinct methods: one turns rainbow Tur\'an numbers into geometric bounds, while the other assigns linear sums to tuples of grid points and exploits the affine relations forced by equal sums.

\bigskip
\noindent\textbf{Notation. }Throughout, we only consider simple undirected graphs. For any graph $G=(V,E)$, write $v(G)=|V(G)|$, $e(G)=|E(G)|$, and $c(G)$ for its number of connected components. For each $v\in V(G)$, write $d_G(v)$ for its degree in $G$. An $m\times n$ matrix $M$ with $m\le n$ is said to be a \emph{full spark} matrix if every $m\times m$ submatrix of $M$ is invertible. An \emph{affine subspace} of $\mathbb R^D$ is a translate $a+V$ of a linear subspace $V\subseteq \mathbb R^D$ for some $a\in \mathbb R^D$, and its dimension is $\dim V$. For a non-empty set $S\subseteq \mathbb R^D$, its \emph{affine hull}, denoted by $\aff(S)$, is the smallest affine subspace containing all the elements of $S$. Equivalently,
$$\aff(S)=\bigg\{\sum_{j=1}^m\lambda_ja_j:a_i\in S,\forall1\le i\le m,\sum_{j=1}^m\lambda_j=1\bigg\}.$$

\bigskip
\noindent\textbf{Hyperplane construction. }
We first construct a random sparsification of an explicit family of hyperplanes to obtain a bound for any fixed \(d,q\) and sufficiently large \(p\).

\begin{theorem}
\label{thm:planar-main}
For any fixed integers $q>d\ge 2$,
$$\HD_d(p,q)\ge p^{\frac{q+d-2}{q-1}+o(1)}$$
as $p\rightarrow\infty$.
\end{theorem}

This result significantly improves the estimate $\HD_2(p,q)\ge p^{1+1/(4q-7)-o(1)}$ in \cite{KellerSmorodinskyLower} to $\HD_2(p,q)\ge p^{1+1/(q-1)-o(1)}$ for $q\ge 3$, and also raises the bound $\HD_2(p,3)\ge p^{5/4-o(1)}$ to $\HD_2(p,3)\ge p^{3/2-o(1)}$.

The construction starts from a hyperplane family with Vandermonde-type normal vectors,
and considers the $q$-uniform hypergraph whose edges are $q$-tuples of
hyperplanes concurrent at a finite point.
A probabilistic incidence supersaturation estimate, combined with the terminal container
theorem, shows that every subfamily containing no $q$-fold concurrence lies
inside one of comparatively few small containers.
We then randomly sparsify the hyperplane family.
The container bound keeps the largest $q$-independent subfamily small,
while a probabilistic estimate controls the maximum number of
sampled hyperplanes through one point.
Finally, we truncate the sampled hyperplanes to compactify them without changing
their finite concurrence pattern.
The ratio between the number of retained hyperplanes and their maximum depth
then gives the required lower bound for the piercing number.

\bigskip

\noindent\textbf{Translation to general-positioned grid sets. }
Our next constructions strongly rely on the structure of sufficiently general-positioned grid sets. More precisely,
for integers $D\geq d$ and $q\geq d+1$, we say that a finite set
$S\subseteq\mathbb R^D$ is in $(d-1,q)$-general position if no $q$
points of $S$ lie in an affine subspace of dimension at most $d-1$.
Let $A_{D,d-1,q}(n)$
denote the maximum size of a set
$S\subseteq[n]^D$
in $(d-1,q)$-general position.
We establish the following ingenious connection between HD numbers and the value of $A_{D,d-1,q}(n)$.

\begin{theorem}
\label{thm:HD&grid}
Fix integers $d\geq 2$, $q\geq d+1$, and real number $0<\rho<1$. 
Suppose that for every sufficiently large integer $D$, 
there exists constant $C_D$ 
such that $A_{D,d-1,q}(n)\leq C_Dn^{\rho D}$, then
$$\text{HD}_d(p,q)\geq p^{1/\rho-o(1)}.$$
\end{theorem}

The proof combines a dimension-preserving projection with point--hyperplane duality.
After projecting $[n]^D$ into $\mathbb R^d$, every affine hyperplane contains at most $n^{d-1}$ projected grid points.
If $p>A_{D,d-1,q}(n)$, then every $p$ projected points contain $q$ points on a common hyperplane.
Dualizing gives a family with the $(p,q)$-property and piercing number of order $n^{D-d+1}$.
Choosing $n$ appropriately and then taking $D$ sufficiently large yields the exponent $1/\rho$.

It remains to obtain upper bounds for $A_{D,d-1,q}(n)$.
To this end, we develop the following two methods.

\bigskip
\noindent\textbf{Graph-based bounds. }
We build up relations between the upper bound of $A_{D,d-1,q}$ and a rainbow Tur\'an bound. Let $\ex^*(n,H)$ denote the largest possible number of edges in a graph on $n$ vertices with a proper coloring and without a rainbow copy of $H$.

\begin{theorem}\label{thm:grid&rbturan}
Fix $d\geq2$ and $q\geq d+1$. Let $H$ be a fixed connected bipartite graph with $e(H)\geq q$ and $v(H)\leq d+1$, then 
$$A_{D,d-1,q}(n)\le\frac{\ex^*(2(e(H)n+1)^D,H)}{(e(H)n+1)^D}.$$ Consequently, if $\mathrm{ex}^*(n,H)=O(n^{1+\rho})$ for some $0<\rho<1$, then
$A_{D,d-1,q}(n)=O_{D,H}(n^{\rho D}),$
and furthermore,
$$\mathrm{HD}_d(p,q)\geq p^{1/\rho-o(1)}.$$
\end{theorem}
Starting from any given $(d-1,q)$-general positioned grid set $S$, the proof uses two copies of sufficiently large grids as the host vertex set, forms an edge whenever the difference of the two vertices falls in $S$, and colors the edge according to the exact element. 
Appearance of a rainbow copy of $H$ produces many affine relations among the colors used by the copy, and forces them to lie in an affine subspace of dimension at most $d-1$, contradicting the $(d-1,q)$-general position condition. Hence the corresponding rainbow Tur\'an bound is also an upper bound of $A_{D,d-1,q}$.

Therefore we aim at optimization of the rainbow Tur\'an exponent under the given parameter restriction. To this end, probably of independent interest, we explore the transition of rainbow Tur\'an behavior under suspension operations.
Here, for a bipartite graph $H=(V_1\cup V_2,E)$, the \emph{suspension} of $H$, denoted by $\Sigma H$, is obtained by adding two adjacent vertices, one to each part, and joining each new vertex to every original vertex in the opposite part. Formally,
$$V(\Sigma H)=(V_1\cup\{x\},V_2\cup\{y\}),$$
$$E(\Sigma H)=E(H)\cup\{xu:u\in V_2\}\cup\{vy:v\in V_1\}\cup\{xy\}.$$
Write $\Sigma^jH$ for the $j$-time suspension of $H$, iteratively defined by $\Sigma^jH=\Sigma(\Sigma^{j-1}H)$. 

\begin{theorem}\label{thm:rbturan&suspension}
 Let \(F\) be a fixed connected bipartite and suppose \(\operatorname{ex}^{*}(n,F)=O(n^{1+\alpha})\) for some \(0\leq\alpha<1\). Then
\[
\operatorname{ex}^{*}(n,\Sigma F)=O(n^{1+\frac{1}{2-\alpha}}).
\]   
\end{theorem}

Starting from the theta graph $\Theta_{k,t}$, whose rainbow Tur\'an behavior was introduced by Janzer~\cite{Janzer2023}, we properly add pendant vertices and apply multiple suspension operations to produce a refined feasible bipartite. Plugging it in Theorem~\ref{thm:grid&rbturan} gives us the optimized exponent $\Gamma(d,q)$, formulated as follows. 

When $q>\lfloor(\frac{d+1}{2})^2\rfloor$, simply write $\Gamma(d,q)=1$. For integer $0\leq s\leq\lfloor\frac{d-1}{2}\rfloor$, write $q_s=(s+1)(d-s)$.
For $d+1\leq q\leq\lfloor(\frac{d+1}{2})^2\rfloor$, there exists a unique integer $s=s(q)$ such that 
$q\in(q_s,q_{s+1}]$. 
Put $t=q-q_s+1$, $k=1+\lfloor\frac{d-2s-1}{t}\rfloor$, 
and define
$$\Gamma(d,q)=\frac{(s+1)k-s}{sk-s+1}.$$

\begin{theorem}\label{thm:graph-main}
For every fixed $d\geq 3$ and $d+1\leq q\leq\lfloor(\frac{d+1}{2})^2\rfloor$,
$$\mathrm{HD}_d(p,q)\geq p^{\Gamma(d,q)-o(1)}.$$
\end{theorem}

In particular, $\Gamma(d,d+1)=\floor{\frac{d+1}2}$, so $\HD_d(p,d+1)\ge p^{\floor{\frac{d+1}2}-o(1)}$. When $d\ge 5$, this significantly improves the estimate $\HD_d(p,d+1)\ge p^{2-o(1)}$ in \cite{KellerSmorodinskyUnions2026}.

\bigskip
\noindent\textbf{Linear reduction of tuples. }
We also establish the following upper bounds of $A_{D,d-1,q}(n)$ by linearly mapping tuples in any $(d-1,q)$-general positioned set to shorter tuples, and uniformly bounding the sizes of preimages at each image tuple.

\begin{theorem}\label{thm:fingerprint-main}
Fix integers \(d\geq2\), \(q\geq d+1\), \(t\geq2\), and \(1\leq s\le h\). If \(th\leq q\) and \(s(t-1)\geq q-d\), then, for every integer \(D\geq d\),
\[
A_{D,d-1,q}(n)=O_{D,d,q,t,h,s}(n^{\frac{sD}{h}}).
\]
Consequently,
\[
\HD_{d}(p,q)\geq p^{\frac{h}{s}-o(1)}.
\]    
\end{theorem}

To derive the bound, given any $(d-1,q)$-general positioned set $S\subseteq [n]^D$, and any $h$-tuple consisting of distinct elements of $S$, we map the tuple to $(\mathbb Z^D)^s$ using an $(s\times h)$ full spark matrix $\Lambda$. Any point in $(\mathbb Z^D)^s$ with too large preimage would result in many affine relations, and ultimately contradicts the $(d-1,q)$-general position condition. Estimating both the amount of possible images and the size of the preimages gives us the desired bounds.

Under the restrictions of the parameters, we deduce the optimal exponent $\Pi(d,q)$, formulated as 
$$\Pi(d,q)=1+\max_{s\in\{1,2,\cdots,d-1\}}\frac 1s\floor{\frac{d-s}{1+\ceil{(q-d)/s}}}.$$

\begin{cor}\label{cor:fingerprint-exp}
For any fixed integers $q\geq d+1\geq 3$, 
$$\HD_d(p,q)\geq p^{\Pi(d,q)-o(1)}$$
\end{cor}

Finally, we compare the exponents $\Gamma(d,q)$ and $\Pi(d,q)$ and discuss the cases in which they are identical. We show that the numerical coincidences reflect a deeper structural correspondence.

\end{section}

\begin{section}{Hyperplane Construction via Containers}
In this section, we prove Theorem \ref{thm:planar-main} by constructing a hyperplane family. We define the base hyperplane family in Section \ref{sec:planar-def}, on which we build an auxiliary hypergraph whose edges are the $q$-concurrences of the hyperplanes, and discuss its properties. We then use the probability-distribution type container introduced by \cite{CamposSamotij2026} to construct a container family for the independent sets in the hypergraph in Section \ref{sec:planar-container}. By randomly sparsifying and truncating the base family in Section \ref{sec:planar-sampling}, we obtain the final witness, where the containers control the absences of $q$-concurrences and provide the $(p,q)$-properties.

\subsection{The auxiliary hypergraph and its properties}\label{sec:planar-def}
Let $q>d\ge 2$ be fixed integers, while $m$ is a prime tending to infinity. Set $n=\ceil{m(m\log m)^{\frac{q-1}{d-1}}}$. 
For each $t\in\mathbb F_m$, define $v_t=(\overline t,\overline{t^2},\cdots,\overline{t^{d-1}})\in \mathbb F_m^{d-1}$, where the bars denote the residues modulo $m$.
For each $a\in[m]$, $b\in[n]$, set $H_{a,b}=\{(\bm x,y):\bm x\in \mathbb R^{d-1},y\in \mathbb R,y=v_a\cdot \bm x+b\}$, and set 
$\cL=\{H_{a,b}:a\in[m],b\in[n]\}.$
Thus \(\abs{\cL}=mn\).  

For any $d$ hyperplanes $H_{a_1,b_1},\cdots,H_{a_d,b_d}$ with $a_1,\cdots,a_d$ pairwise distinct, the matrix formed by their normal vectors $(v_a~~-1)$ is invertible, since its determinant equals $(-1)^d\prod\limits_{i<j}(a_j-a_i)\ne0$ modulo $m$. It follows that any $d$ hyperplanes in $\cL$ with distinct directions intersect at exactly one point, and that any $q$ hyperplanes intersect at at most one point. We then define the auxiliary $q$-uniform hypergraph $\cH$ on $\cL$, where $q$ hyperplanes form an edge if and only if they intersect at one point. 

To explore the intersection properties of $\cL$, we consider its incidence with a local grid set. Setting $\cP=[\floor{\frac nm}]^{d-1}\times [dn]\subseteq \mathbb Z^d$, we have the following elementary results.

\begin{lemma}\label{lem:planar-elem-H}
\begin{enumerate}[(a)]
\item For each $H\in \cL$, $|H\cap \cP|=\floor{\frac nm}^{d-1}$.
\item For each $j\ge1$ and pairwise distinct hyperplanes $H_1,\cdots,H_j\in \cL$, $|(\bigcap\limits_{i=1}^jH_i)\cap\cP|\le(\frac{n}{m})^{\max\{d-j,0\}}$.
\end{enumerate}
\end{lemma}

\begin{proof}
To see (a), note that for each fixed $H=H_{a,b}\in\cL$ and every $\bm x\in[\floor{\frac nm}]^{d-1}$, we have $0\le v_a\cdot \bm x+b\le (d-1)m\floor{\frac nm}+n\le dn$. Therefore, each $\bm x\in [\floor{\frac nm}]^{d-1}$ uniquely determines a coordinate $y=v_a\cdot \bm x+b\in[dn]$, and then an element of $H\cap \cP$. Thus $|H\cap \cP|=\floor{\frac nm}^{d-1}$.

To see (b), first note that when $j\ge d$, the hyperplanes $H_1,\cdots,H_j$ intersect in at most one point, and the inequality clearly holds. When $j<d$, we show that, for any $\lambda_j,\lambda_{j+1},\cdots,\lambda_{d-1}\in[\floor{\frac nm}]$, the hyperplanes $P_i:x_i=\lambda_i$ ($i\in[j,d-1]$) and $H_1,\cdots,H_j$ intersect at most one point. Indeed, suppose $H_i=H_{a_i,b_i}$ for each $i\in [j]$, if $a_{i_1}=a_{i_2}$ holds for some $i_1<i_2$, then $H_{i_1}$ and $H_{i_2}$ are parallel and they do not intersect. Otherwise $a_1,\cdots,a_j$ are pairwise distinct, then the matrix formed by the normal vectors of $H_1,\cdots,H_j,P_j,\cdots,P_{d-1}$ has determinant $(-1)^d\prod\limits_{1\le i_1<i_2\le j}(a_{i_2}-a_{i_1})\ne0$ modulo $m$, thus it is invertible, and the $d$ hyperplanes intersect at most one point. It follows that, when the coordinates $x_j,\cdots,x_{d-1}$ are fixed, there is at most one choice of $(\bigcap\limits_{i=1}^jH_i)\cap\cP$. Since the coordinates $(x_j,\cdots,x_{d-1})$ have $\floor{\frac nm}^{d-j}$ choices, we conclude that $|(\bigcap\limits_{i=1}^jH_i)\cap\cP|\le\floor{\frac nm}^{d-j}$, as desired.
\end{proof}

Using Lemma \ref{lem:planar-elem-H} we establish the following regularity property of $\cH$. 
For conciseness, set $A=(m^{d-1}n^{q-d})^{\frac{1}{q-1}}$, then $A\le\frac{n}{m\log m}$ and $A=(1+o(1))\frac{n}{m\log m}$ as $m$ tends to infinity by definition.

\begin{lemma}\label{lem:planar-supersaturation}
There exists a constant $N=N(d,q)>0$ such that every
\(\cJ\subseteq\cL\) with \(\abs \cJ=M\geq 4dqn\) attains a probability distribution on $E(\cH[\cJ])$ whose resulting edge $R$ satisfies 
$$\Pr(S\subseteq R)\le \frac{NA^{j-1}}{M^j}$$
for every $S\subseteq\cJ$ with $|S|=j\in[2,q]$.

\end{lemma}

\begin{proof}

For a grid point $x\in\cP$, let $d_{\cJ}(x)=|\{H\in \cJ:x\in H\}|$. Double-counting on the point-hyperplane incidences gives 
$$\sum_{x\in\cP}d_{\cJ}(x)=\sum_{H\in \cJ}|H\cap\cP|=M\floor{\frac nm}^{d-1}.$$
On the left-hand-side, all grid points $y$ with $d_{\cJ}(y)<\frac{M}{2dn}$ contribute at most $\frac{M}{2dn}|\cP|\le\frac M2\floor{\frac nm}^{d-1}$.
Consequently, with $\cP':=\{x\in\cP: d_\cJ(x)\geq \frac{M}{2dn}\}$,
\begin{align}\label{eq:planar-total-weight}
\sum_{x\in\cP'}d_{\cJ}(x)\ge\frac{M}{2}\floor{\frac nm}^{d-1}\ge\frac M2(\frac{n}{2m})^{d-1}.
\end{align}

Set $W=\sum\limits_{x\in\cP'}d_{\cJ}(x)$. We define a probability distribution on $E(\cH[\cJ])$ as follows. For each $x\in\cP'$, choose $x$ with probability $\frac{d_{\cJ}(x)}{W}$, conditional on which choose uniformly a $q$-subset of the hyperplanes in $\cJ$ containing $x$. The resulting set is indeed an edge of $\cH[\cJ]$, and the probabilities of all choices sum to $1$, so the procedure determines a probability distribution on $E(\cH[\cJ])$. 

Arbitrarily fix $S\subseteq \cJ$ with $|S|=j\in [2,q]$. For each $y\in (\bigcap\limits_{H\in S}H)\cap\cP'$, the probability that $y$ is chosen in the above procedure is $\frac{d_{\cJ}(y)}{W}$. Conditional on $y$ being chosen, there are ${{d_{\cJ}(y)}\choose{q}}$ choices of $R$, and ${{d_{\cJ}(y)-j}\choose{q-j}}$ of them contain $S$ as a subset. Hence
$$\Pr(S\subseteq R|y~\text{is chosen})=\frac{{{d_{\cJ}(y)-j}\choose{q-j}}}{{{d_{\cJ}(y)}\choose{q}}}=\prod_{i=0}^{j-1}\frac{q-i}{d_{\cJ}(y)-i}\le(\frac{2q}{d_{\cJ}(y)})^j$$
using $d_{\cJ}(y)\ge \frac{M}{2dn}\ge 2q>2i$ for each $i\in[0,j-1]$.

Eventually, by Lemma \ref{lem:planar-elem-H}(b), $|(\bigcap\limits_{H\in S}H)\cap\cP'|\le(\frac{n}{m})^{\max\{d-j,0\}}$, thus
\begin{align*}
\Pr(S\subseteq R)&\le \sum_{y\in (\bigcap\limits_{H\in S}H)\cap\cP'}\frac{d_{\cJ}(y)}{W}(\frac{2q}{d_{\cJ}(y)})^j\\
&\le (\frac{n}{m})^{\max\{d-j,0\}}(2q)^j(\frac{M}{2dn})^{1-j}\frac{2}{M}(\frac{2m}{n})^{d-1}\\
&\le(2^{2j+d}q^jd^j)M^{-j}n^{j-1}(\frac{m}{n})^{\min\{d-1,j-1\}}\\
&\le(2^{2j+d}q^jd^j)M^{-j}n^{j-1}(\frac{m}{n})^{\frac{(j-1)(d-1)}{q-1}}\\
&=(2^{2j+d}q^jd^j)\frac{A^{j-1}}{M^j},
\end{align*}
where the second inequality comes from \eqref{eq:planar-total-weight}, and the fourth inequality comes from $\frac mn<1$ and $\min\{d-1,j-1\}\ge \frac{(d-1)(j-1)}{q-1}$. Setting $N=2^{2q+d}q^qd^q$ completes the proof.
\end{proof}

\end{section}

\subsection{Terminal containers}\label{sec:planar-container}
We use the probabilistic distribution form of Theorem~D of Campos and Samotij~\cite{CamposSamotij2026}, with its
parameter renamed \(\theta\).

\begin{theorem}[Campos--Samotij~\cite{CamposSamotij2026}, weighted form of Theorem~D]
\label{thm:planar-campos-samotij}
Let \(G\) be an \(r\)-uniform hypergraph on a finite vertex set \(V\), and
let \(0<\theta\leq1/(8r^2)\).  There is a family
\(\mathcal S\subseteq2^V\) with maps
\(g:\cI(G)\rightarrow\mathcal S\) and \(f:\mathcal S\rightarrow2^V\)
satisfying the following properties.
\begin{enumerate}[label=\textup{(\roman*)}]
\item
For every \(I\in\cI(G)\),
\(
  g(I)\subseteq I\subseteq f(g(I)).
\)
\item
Every \(S\in\mathcal S\) satisfies
\(
  \abs S\leq8r^2\theta\abs V.
\)
\item
Suppose $C=f(S)\ne \emptyset$ for some $S\in \mathcal S$, then for every probability distribution on $E(G[C])$ with random edge $R$, there is a set $U\subseteq C$ with $2\le |U|\le r$ satisfying $\Pr(U\subseteq R)\ge\frac{\theta^{|U|-1}}{|C|}$.
\end{enumerate}
\end{theorem}

We first construct a one-step container family for each large subset of $\cL$.

\begin{lemma}\label{lem:planar-one-step-container}
For all sufficiently large $m$, and any $\cJ\subseteq \cL$ satisfying $|\cJ|=M\ge8dqn$, there is a family $\cC(\cJ)\subseteq 2^{\cJ}$, such that
\begin{enumerate}[label=\textup{(\roman*)}]
\item Any independent set in $\cJ$ is contained in some $C\in \cC(\cJ)$.
\item Every $C\in\cC(\cJ)$ has size less than $\frac M2$.
\item $|\cC(\cJ)|<M^{s+1}$, where $s=\ceil{32q^2NA}$.
\end{enumerate}
\end{lemma}

\begin{proof}
Set $\theta=\frac{4NA}{M}$, then $0<\theta\le\frac{2NA}{dqn}\le \frac{2N}{qm\log m}<\frac{1}{8q^2}$ for sufficiently large $m$. Therefore we can apply Theorem \ref{thm:planar-campos-samotij} with $G=\cH[\cJ]$ and $\theta$ chosen above to output a fingerprint family $\mathcal S(\cJ)$ and a container family $\cC(\cJ)=f(\mathcal S(\cJ))$, such that any independent set is contained in some $C\in\cC(\cJ)$.

We first prove that each $C\in \cC(\cJ)$ satisfies $|C|< \frac M2$. Assume to the contrary that $|C|\ge \frac M2\ge 4dqn$, we apply Lemma \ref{lem:planar-supersaturation} to $E(\cH[C])$ to obtain a probability distribution whose random edge $R$ satisfies
$$\Pr(U\subseteq R)\le \frac{NA^{j-1}}{|C|^j}\le(\frac{2NA}{M})^{j-1}\frac{1}{|C|}<\frac{\theta^{j-1}}{|C|}$$
for each $U\subseteq C$ with $|U|=j\in[2,q]$, which contradicts Theorem \ref{thm:planar-campos-samotij}.

It remains to show that $|\mathcal C(\cJ)|<M^{s+1}$. For each $S\in\mathcal S(\mathcal J)$, $|S|\le 32q^2NA\le s$ by Theorem~\ref{thm:planar-campos-samotij}(\romannumeral2). Therefore
$$|\cC(\cJ)|\le|\mathcal S(\cJ)|\le\sum_{i=0}^s{M\choose i}<M^{s+1}.$$
\end{proof}

We then iteratively apply Lemma~\ref{lem:planar-one-step-container} to construct a global container family.

\begin{lemma}
\label{lem:planar-terminal-containers}
There is a constant \(B=B(d,q)>0\) such that, for all sufficiently large
\(m\), there is a family \(\cC\) of subsets of \(\cL\) satisfying the
following properties.
\begin{enumerate}[label=\textup{(\roman*)}]
\item Every independent set of \(\cH\) is contained in some
\(C\in\cC\).
\item Every \(C\in\cC\) has size less than \(8dqn\).
\item $\log\abs{\cC}\leq BA(\log m)^2$.
\end{enumerate}
\end{lemma}

\begin{proof}
We run the following algorithm to obtain the container family $\cC$.
\begin{algorithm}[!h]
\caption{Construction of Global Containers}
\label{alg:container}
\renewcommand{\algorithmicrequire}{\textbf{Input:}}
\renewcommand{\algorithmicensure}{\textbf{Output:}}
\begin{algorithmic}[1]
\REQUIRE $\mathcal H$
\ENSURE $\cC$
\STATE $i\gets0$, $\cC_0\gets \{V(\mathcal H)\}$
\WHILE{$\exists \cJ\in \cC_i$,~$|\cJ|\ge8dqn$}
\STATE Obtain $\cC(\cJ)$ by Lemma~\ref{lem:planar-one-step-container}
\STATE $\cC_{i+1}\gets(\cC_i\backslash\{\cJ\})\cup\cC(\cJ)$, $i\gets i+1$
\ENDWHILE
\RETURN $\cC_i$
\end{algorithmic}
\end{algorithm}

To see that the algorithm terminates in finitely many steps, set $\tau=\ceil{(s+1)\log_2 mn}$, and define 
$$w_i=\sum_{C\in \cC_i}|C|^{\tau}.$$
Then $$w_i-w_{i+1}\ge|\cJ|^{\tau}-\sum_{C\in \cC(\cJ)}|C|^{\tau}>|\cJ|^{\tau}-|\cJ|^{s+1}(\frac{|\cJ|}{2})^{\tau}> 0,$$
where $\cJ$ is the subset chosen at the $i$-th step.
Since every $w_i$ is a positive integer, the algorithm terminates in $w_0$ steps.

The first desired property comes from Lemma~\ref{lem:planar-one-step-container}(\romannumeral1), and the independent sets are inherited at each step. The second directly comes from the loop condition. To prove the third property, suppose the algorithm terminates at the $t$-th step, and note that every set in $\cC$ except the empty container, if one exists, contributes at least $1$ to $w_t$. Therefore
$$|\cC|\le w_t+1<w_0+1=|V(\cH)|^{\tau}+1,$$
and 
$$\log|\cC|<\tau\log|V(\cH)|+1\le2s(\log(mn)^2)<64q^2N A(4q\log m)^2.$$
Pick $B=2^{10}q^4N$ and we complete the proof.
\end{proof}

\begin{subsection}{Sparse Sampling}\label{sec:planar-sampling}

We proceed by randomly sparsifying $\cL$ with properly chosen parameters. The terminal container family constructed above allows us to control all subfamilies containing no $q$ concurrent hyperplanes simultaneously. Meanwhile, an elementary union-bound argument controls the maximum number of sampled hyperplanes through any point. Together, these estimates yield a hyperplane family with the desired intersection property and large transversal number.

Let $X$ be a random subfamily of $\cL$, and each hyperplane in $\cL$ appears in $X$ independently with probability $\frac {\log m}{m}$. Let $\alpha_q(X)$ be the maximal size of a subfamily of $X$ in which no $q$ hyperplanes share a common point. Equivalently, $\alpha_q(X)$ is the maximal size of an independent set in $\cH[X]$. Define the maximal depth $\Delta(X)$ as the maximal number of hyperplanes in $X$ through a common point. Formally, $\Delta(X)=\max\limits_{x\in\mathbb R^d}|\{\ell\in X:x\in\ell\}|$.

\begin{theorem}\label{thm:planar-sampling}
There exist constants $D=D(d,q)$, $K=K(d,q)$ such that, with probability at least $\frac 12$, $X$ satisfies all of the following conditions.
\begin{enumerate}[label=\textup{(\roman*)}]
\item $|X|\ge \frac{n\log m}{2}$.
\item $\alpha_q(X)<p_m:=\ceil{DA\log^2m}$.
\item $\Delta(X)<K\log m$.
\end{enumerate}
\end{theorem}

We will use the following elementary binomial estimates.

\begin{lemma}
\label{lem:binomial-tails}
For any integer $N$ and any probability $0<\rho\leq 1$, let \(Y\sim\operatorname{Bin}(N,\rho)\) and let
\(\mu=N\rho\). 
\begin{enumerate}[(a)]
\item
  \(\Pr(Y<\mu/2)\leq \exp(-\mu/8)\).
\item
  If \(t\) is an integer with \(t\geq\mu\), then
  \(\Pr(Y\geq t)\leq
    \left(\frac{e\mu}{t}\right)^t.\)
\end{enumerate}
\end{lemma}

\begin{proof}
For the first estimate, Markov's inequality gives
\[
  \Pr(Y<\mu/2)
  \leq
  2^{\mu/2}\E[2^{-Y}]
  =
  2^{\mu/2}(1-\rho/2)^N
  \leq
  \exp\left(-\frac{1-\log 2}{2}\mu\right)
  \leq
  \exp(-\mu/8).
\]
If \(t>N\), the second estimate is immediate because the event is empty.
Otherwise, the event \(Y\geq t\) implies that some
\(t\)-subset of the \(N\) trials consists entirely of successes.  A union bound over all possible choices of the $t$-subset and the elementary bound \(\binom Nt\leq(eN/t)^t\) yield
\[
  \Pr(Y\geq t)
  \leq
  \binom Nt\rho^t
  \leq
  \left(\frac{eN\rho}{t}\right)^t.
\]
\end{proof}

\begin{proof}[Proof of Theorem~\ref{thm:planar-sampling}]
By Lemma~\ref{lem:binomial-tails}(a) with $\mu=n\log m$, we have 
\begin{align}\label{eq:planar-X-large}
\Pr(|X|<\frac{n\log m}{2})\le\exp\left(-\frac{n\log m}{8}\right).
\end{align}
Let $\cC$ be the terminal container family obtained by Lemma~\ref{lem:planar-terminal-containers}. For each $C\in\cC$, 
$$\mathbb E(|X\cap C|)=\frac {\log m}m|C|\le \frac{8dqn\log m}{m}\le20dqA\log^2m$$
by Lemma~\ref{lem:planar-terminal-containers}(\romannumeral2).
Pick $D=200(B+1)dq$, where $B$ is the constant in Lemma~\ref{lem:planar-terminal-containers}(\romannumeral3), 
then by Lemma~\ref{lem:binomial-tails}(b),
$$\Pr(|X\cap C|\ge p_m)<(\frac{20edq}{D})^{p_m}<\exp(-p_m).$$
If $\alpha_q(X)\ge p_m$, that is, $\cH[X]$ contains an independent set of size at least $p_m$, then by Lemma~\ref{lem:planar-terminal-containers}(\romannumeral1), there exists some $C\in\cC$ such that $|X\cap C|\ge p_m$.
Therefore by Lemma~\ref{lem:planar-terminal-containers}(\romannumeral3) and a union bound,
\begin{align}\label{eq:planar-alphaX-large}
\Pr(\alpha_q(X)\ge p_m)\le \exp(BA\log^2m-p_m)<\exp(-A\log^2m).   
\end{align}
Now set $\cP''$ to be the set of points which are contained by at least $d$ hyperplanes in $\cL$. Since each point in $\cP''$ is determined by $d$ hyperplanes in $\cL$, we have $|\cP''|\le|\cL|^d=m^dn^d$. For $x\in\cP''$, let $S_x$ denote the number of hyperplanes in $X$ passing through $x$. Since each point $x\in \cP''$ is contained by at most $m$ hyperplanes in $\cL$, $\mathbb ES_x\le \log m$. Thus by setting $K_q=5e^2dq$ and using Lemma~\ref{lem:binomial-tails}(b),
$$\Pr(S_x\ge K_q\log m)\le \exp(-5dq\log m)=m^{-5dq}.$$
Take a union bound over all elements in $\cP''$,
\begin{align}\label{eq:planar-deltaX-large}
\Pr(\Delta(X)\ge K_q\log m)\le m^dn^dm^{-5dq}<\frac 1m.  \end{align}
Combining \eqref{eq:planar-X-large}\eqref{eq:planar-alphaX-large}\eqref{eq:planar-deltaX-large} and a final union bound yields the conclusion.
\end{proof}

Pick a subfamily $X$ satisfying the three conditions listed above. To finally obtain a compact family, we choose a sufficiently large closed ball $\mathcal D$, such that every finite incidence among $X$ intersects $\mathcal D$, and truncate $X$ with $\mathcal D$ to obtain a family $X'$ of compact convex sets. That is $X':=\{H\cap \mathcal D:H\in X\}$.

\begin{cor}\label{cor:planar-truncate}
The truncated set $X'$ satisfies the following.   \begin{enumerate}[(a)]
\item $X'$ has $(p_m,q)$-property.
\item $\tau(X')>\frac{n}{2K}.$
\end{enumerate} 
\end{cor}

\begin{proof}
By Theorem~\ref{thm:planar-sampling}(\romannumeral2), among every $p_m$ hyperplanes in $X$, at least $q$ hyperplanes are concurrent. By the choice of $\mathcal D$, the corresponding $q$ truncated hyperplanes are also concurrent. Thus $X'$ has $(p_m,q)$-property.

Additionally, each point is contained in less than $K\log m$ hyperplanes in $X$, so $\tau(X)\ge\frac{|X|}{K\log m}>\frac{n}{2K}$. Since each transversal of $X'$ is also a transversal of $X$, we conclude that $\tau(X')>\frac{n}{2K}$.
\end{proof}

\begin{proof}[Proof of Theorem~\ref{thm:planar-main}]
For sufficiently large $p$, we approximate $p$ with some $p_m$. Let $k$ be the largest integer such that $Dk^{\frac{q-1}{d-1}}(\log k)^{2+\frac{q-d}{d-1}}<p$, and let $m$ be the maximum prime in $[1,k]$, then $Dk^{\frac{q-1}{d-1}}(\log k)^{2+\frac{q-d}{d-1}}\ge\frac p2$ and $m\ge \frac k2$. By definition of $A$, we have $p_m<p$ and $p_m=\ceil{DA\log^2m}>2^{-q}Dk^{\frac{q-1}{d-1}}(\log k)^{2+\frac{q-d}{d-1}}\ge 2^{-q-1}p$.

Let $X'$ be defined as above. Since $X'$ has $(p_m,q)$-property, it also has $(p,q)$-property. Moreover, as $p$ tends to infinity, combining $p<2^{q+1}p_m=m^{\frac{q-1}{d-1}+o(1)}$ and $\tau(X')>\frac{n}{2K}=m^{\frac{q+d-2}{d-1}+o(1)}$ we have $\tau(X')\ge p^{\frac{q+d-2}{q-1}+o(1)}$. Therefore, as a compact convex family, $X'$ witnesses the bound $\HD(p,q)\ge p^{\frac{q+d-2}{q-1}+o(1)}$, as desired.
\end{proof}

\end{subsection}

\begin{section}{Hadwiger-Debrunner Numbers and General-positioned Grid Sets}

In this section, we relate lower bounds for Hadwiger--Debrunner numbers
to upper bounds on the size of general-positioned grid sets.
A suitable linear projection preserves the relevant affine dependencies
while controlling hyperplane incidences, and point--hyperplane duality
then converts the projected grid into a family of compact convex sets.

We use the following variant of the point--hyperplane duality lemma of Keller and Smorodinsky~\cite{KellerSmorodinskyUnions2026} to turn a point configuration into a family of compact convex sets.

\begin{lemma}
\label{lem:point-configuration-transfer}
Let \(p\geq q\geq d+1\), and let \(P\subseteq\R^{d}\) consist of \(N\geq p\) points. Suppose every \(p\)-point subset of \(P\) contains \(q\) points on a common affine hyperplane, while every affine hyperplane contains at most \(M\) points of \(P\). Then
\[
\HD_{d}(p,q)\geq\frac{N}{M}.
\]
\end{lemma}

\begin{proof}
For each $P_1\subseteq P$ and $|P_1|=p$, choose a hyperplane $H(P_1)$ containing at least $q$ of its elements, and denote the direction space of $H(P_1)$ by $U(P_1)$, then $\dim U(P_1)=d-1$. Then choose a unit vector $v$ such that $\bm v\notin \bigcup\limits_{P_1\subseteq P, |P_1|=p}U(P_1)$. Such a vector $\bm v$ must exist, since we only exclude finitely many proper subspaces of $\mathbb R^d$. Choose a new orthogonal basis of $\mathbb R^d$ using $v$ as the $d$-th basis vector, then every witnessing hyperplane defined above is not parallel to the $x_d$-axis.

We then introduce the following point-hyperplane duality function. For $a\in \mathbb R^d$, define $H(a)$ to be the hyperplane $x_d=\sum\limits_{j=1}^{d-1}a_jx_j-a_d$. For each hyperplane $H:\sum\limits_{i=1}^dh_ix_i=c$ not parallel to $x_d$-axis, which means $h_d\ne 0$, define $P(H)=(-\frac{h_1}{h_d},-\frac{h_2}{h_d},\cdots,-\frac{h_{d-1}}{h_d},-\frac{c}{h_d})$. One easily checks that $H(a)$ is not parallel to the $x_d$ axis, that $P(H(a))=a$, and that $a\in H$ if and only if $P(H)\in H(a)$, for any point $a$ and hyperplane $H$.

Set $\cH=\{H(a):a\in P\}$. For any $p$ hyperplanes $H(a_1),H(a_2),\cdots H(a_p)\in\cH$, a witnessing hyperplane $H_0$ contains $q$ of $a_1,a_2,\cdots,a_p$, thus $P(H_0)$, whose existence is ensured by the choice of basis, is contained by $q$ of $H(a_1),H(a_2),\cdots H(a_p)$. Hence $\cH$ has $(p,q)$-property. Moreover, if any point $p$ is contained by $M+1$ hyperplanes $H(a_1),H(a_2),\cdots H(a_{M+1})$, then $a_1,a_2,\cdots,a_{M+1}\in H(p)$, a contradiction. Therefore any point lies in at most $M$ hyperplanes in $\cH$. To obtain a compact family, we truncate $\mathcal H$ with a sufficiently large closed ball that intersects every incidence of $M$ hyperplanes to obtain a new family $\cH'$, then $\cH'$ also has $(p,q)$-property, and $\tau(\cH')\ge\frac{|\cH|}{M}=\frac NM$. We conclude that $\HD_{d}(p,q)\ge \frac NM$. 
\end{proof}

The following lemma provides the point set with desired point-hyperplane incidences. As below, for a linear mapping $A$, we also write $A$ for its matrix.

\begin{lemma}\label{lem:HD&A-project-image}
For all positive integers $D>d$ and $n$, there exists an injective linear mapping $\pi:[n]^D\rightarrow\mathbb R^d$ whose image $P=\pi([n]^D)$ has the property that every hyperplane in $\mathbb R^d$ contains at most $n^{d-1}$ points of $P$.
\end{lemma}
\begin{proof}
We first find a mapping $\pi$ preserving the affinely independent relations among no more than $d+1$ points. Set
$$\cT=\{T\subseteq[n]^D:|T|\le d+1,~\text{elements of}~T~\text{are affinely independent}\}.$$
For each $T=(u_0,u_0+u_1,\cdots,u_0+u_{t})\in\cT$, where $t\le d$, we extend $(u_1,\cdots,u_t)$ to a basis $(u_1,\cdots,u_D)$ of $\mathbb R^D$, and for any $(d\times D)$-matrix $A$, set 
$$p_T(A)=\det\left(Au_1~Au_2~\cdots~Au_d\right).$$
Then $p_T$ is a polynomial on the $(d\times D)$-matrix space. Also note that, if $p_T(A)\ne0$, then the image of $T$ under $A$ is still affinely independent. We also show that $p_T\ne0$. Actually, consider $A_T$ which maps $u_i$ to $e_i$, the $i$-th standard basis vector of $\mathbb R^d$, for each $1\le i\le d$, then $p_T(A_T)=1$. Therefore, set $Q(A)=\prod\limits_{T\in\cT}p_T$, and $Q$ is also non-zero. Pick $\pi$ such that $Q(\pi)\ne0$, we have $p_T(\pi)\ne0$ for each $T$, so $\pi$ preserves the affine independence of all tuples in $\cT$ simultaneously. In particular, $\pi u_1\ne \pi u_2$ for any $u_1\ne u_2$, so $\pi$ is injective.

Set $P=\pi([n]^D)$. Fix a hyperplane $H\subset \mathbb R^d$, set $X=\pi^{-1}(P\cap H)$. Fix any $x_0\in X$ and let $W$ denote the affine subspace spanned by $\{x-x_0:x\in X\}$, with $\dim W:=r$. By the choice of $\pi$, $\dim W\le \dim H=d-1$. Pick a basis $(w_1,\cdots,w_r)$ of $W$, and form a matrix $B=\left(w_1,\cdots,w_r\right)$, then $\rank B=r$, and we assume that the first $r$ rows are linearly independent. 

Let $\varphi:[n]^D\rightarrow[n]^r$ be the projection onto the first $r$ coordinates, so $\varphi=\left(I_r~0\right)$ as a matrix. We show that $\varphi$ is injective on $X$. Indeed, if $\varphi(x)=\varphi(y)$, note that $x-y\in W$, so there exists $v\in \mathbb R^r$ such that $x-y=Bv$, then $\varphi Bv=0$. However by the assumption above $\varphi B$ is invertible, so $v=0$ and $x=y$. Therefore, $|X|\le |[n]^r|=n^r.$
Eventually, $|P\cap H|=|\pi(X)|=|X|\le n^r$. This holds for any hyperplane $H\subset \mathbb R^d$, and we complete the proof.
\end{proof}

\begin{proof}[Proof of Theorem~\ref{thm:HD&grid}]
For any fixed sufficiently large integer $D$, 
and $p$ sufficiently large with respect to $D$, 
set $n$ to be the maximal integer such that $A_{D,d-1,q}(n)<p$, 
then $p\le C_D(n+1)^{\rho D}<2C_Dn^{\rho D}<n^D$. 
Apply Lemma~\ref{lem:HD&A-project-image} to obtain the linear mapping $\pi$, and image $P$. 
We check that $P$ satisfies the conditions in Lemma~\ref{lem:point-configuration-transfer} with properly set parameters. 
For any $p$ points in $P$, their preimages in $[n]^D$ are not in $(d-1,q)$ general position, 
so $q$ of them lie in an affine subspace with dimension $d-1$, 
and the images of these $q$ points lie in a hyperplane of $\mathbb R^d$ by linearity of $\pi$. 
Additionally, every hyperplane contains at most $n^{d-1}$ points of $P$. 
Apply Lemma~\ref{lem:point-configuration-transfer} with $N=n^D,M=n^{d-1}$, 
we have 
$$\HD_d(p,q)\ge n^{D-d+1}>\left(\frac{p}{2C_D}\right)^{\frac{D-d+1}{\rho D}}.$$

To obtain the asymptotic conclusion, fix $\varepsilon>0$, 
and choose $D$ sufficiently large so that the assumed estimates
hold and that $\frac{d-1}{\rho D}<\frac{\varepsilon}{2}$. 
For this fixed $D$ and every sufficiently large $p$, we have
$$\HD_d(p,q)>\left(\frac{p}{2C_D}\right)^{\frac{D-d+1}{\rho D}}>\left(\frac{p}{2C_D}\right)^{\frac{1}{\rho}-\frac{\varepsilon}{2}}>p^{\frac{1}{\rho}-\varepsilon}.$$
Since $\varepsilon$ is chosen arbitrarily, we conclude that $\HD_d(p,q)\ge p^{\frac{1}{\rho}-o(1)}$ as $p\rightarrow\infty$.
\end{proof}

\end{section}

\begin{section}{Rainbow Tur\'an, Suspension and Graph-based Bounds}

In this section, we translate rainbow Tur\'an estimates into upper bounds for general-position grid sets and hence into lower bounds for Hadwiger--Debrunner numbers. In Section \ref{sec:graph-rt-grid}, we first construct, from a grid set in $(d-1,q)$ general position, a properly colored bipartite graph whose edge colors are the points of the grid set. The rational cycle space of a forbidden graph $H$ then yields affine relations among the colors of every rainbow copy of $H$. We next study how rainbow Tur\'an exponents behave under the suspension operation in Section \ref{sec:graph-suspension}. Finally, in Section \ref{sec:graph-quant}, by combining the suspension theorem with Janzer's estimates for theta graphs and a pendant-vertex extension, we obtain the quantitative graph-based bounds.

\begin{subsection}{Rainbow Tur\'an and general positioned grid sets}\label{sec:graph-rt-grid}
Given any $(d-1,q)$-general positioned set $S$ and connected bipartite $H$ with restrictions given in \ref{thm:grid&rbturan}, we construct a host properly colored bipartite graph, where each edge corresponds to an element of $S$, and the existence of any $q$-colored copy of $H$ is forbidden.

We first recall the \emph{rational cycle space} of a bipartite $F$. Choose a bipartition $V(F)=V_L\cup V_R$, and define the incidence matrix $M\in\mathbb Q^{V(F)\times E(F)}$, 
$$M_{u,e}=
\begin{cases}
1,~\text{if}~u\in e~\text{and}~u\in V_L,\\    
-1,~\text{if}~u\in e~\text{and}~u\in V_R,\\  
0,~\text{if}~u\notin e.    
\end{cases}$$
The \emph{rational cycle space} of $F$ is defined by $Z_1(F;\mathbb Q)=\ker M$ over the field of rational numbers. For a cycle $C=u_1w_1\cdots u_kw_ku_{k+1}\subseteq F$, where $u_{k+1}=u_1$, and $u_i\in V_L$, $w_i\in V_R$, we define the cycle vector $\chi_C$ by $(\chi_C)_{u_jw_j}=1$, $(\chi_C)_{u_{j+1}w_j}=-1$, and $(\chi_C)_{e}=0$ for $e\notin E(C)$.

We will use the following results on cycle spaces, see~\cite{Bapat2014}.
\begin{lemma}\label{lem:graph-cycle-space}
\begin{enumerate}[(a)]
\item $\rank(M)=v(F)-c(F)$, $\dim Z_1(F;\mathbb Q)=e(F)-v(F)+c(F)$,
\item Arbitrarily choose a spanning forest $T$ of $F$, and for each $e\in E(F)\backslash E(T)$, let $C(e)$ be the unique cycle in the subgraph $T\cup\{e\}$. Then $\{\chi_{C(e)}:e\in E(F)\backslash E(T)\}$ forms a basis of $Z_1(F;\mathbb Q)$.
\end{enumerate}    
\end{lemma}
\begin{proof}[Proof of Theorem~\ref{thm:grid&rbturan}]
Let \(S\subseteq[n]^D\) be in \((d-1,q)\) general position. Set
\(N=e(H)n+1\) and \(G=(\Z/N\Z)^{D}\), with addition in each coordinate modulo \(N\).
Form a bipartite graph \(\Gamma_{S}\) with two copies \(G_{L},G_{R}\) of \(G\).
For every \(x\in G\) and \(s\in S\), join \(x\in G_{L}\) to
\(x+s\in G_{R}\) and give this edge color \(s\).   The coloring is proper because the colors incident with \(x\) are the distinct elements of \(S\).  

For the sake of contradiction, suppose $\Gamma_S$ contains a rainbow copy of $H$. For each $e\in E(H)$, let $s(e)$ denote the color of $e$. Fix a bipartition $V(H)=V_L\cup V_R$ where $V_L\subseteq G_L$, $V_R\subseteq G_R$. To deduce the affine relations among the colors used by $H$, define a matrix
$A\in \mathbb Z^{(D+1)\times E(H)}$ whose column indexed by $e$ is $(1~~s(e)^{\mathsf T})^{\mathsf T}$.

We show that for each cycle $C=u_1w_1\cdots u_kw_ku_{k+1}\subseteq H$, where $u_{k+1}=u_1$, $u_i\in V_L$, $w_i\in V_R$, the cycle vector $\chi_C$ lies in $\ker A$. It suffices to prove that $\sum\limits_{e\in E(C)}(\chi_C)_e=0$, and that $\sum\limits_{e\in E(C)}(\chi_C)_es(e)=0$. The first assertion follows from 
$$\sum\limits_{e\in E(C)}(\chi_C)_e=\sum\limits_{j=1}^k(\chi_C)_{u_jw_j}+\sum\limits_{j=1}^k(\chi_C)_{u_{j+1}w_j}=k-k=0.$$
To prove the second, write $T(C)=\sum\limits_{e\in E(C)}(\chi_C)_es(e)$, and note that 
$$T(C)=\sum\limits_{e\in E(C)}(\chi_C)_es(e)=
\sum\limits_{j=1}^ks({u_jw_j})-\sum\limits_{j=1}^ks(u_{j+1}w_j)=\sum\limits_{j=1}^k(w_j-u_j)-\sum\limits_{j=1}^k(w_j-u_{j+1})=0$$
holds in $G$. Therefore each coordinate of $T(C)$ is a multiple of $N$. However for each $1\le j\le D$,
$$|T(C)_j|=|(\sum\limits_{e\in E(C)}(\chi_C)_es(e))_j|\le\sum\limits_{e\in E(C)}|s(e)_j|\le ne(H)<N,$$
thus $T(C)_j=0$ for each $j$, and $T(C)=0$.

By Lemma~\ref{lem:graph-cycle-space}(b), a basis of $Z_1(H;\mathbb Q)$ lies in $\ker A$, so $Z_1(H;\mathbb Q)\subseteq \ker A$. By Lemma~\ref{lem:graph-cycle-space}(a), $\dim\ker A\ge \dim Z_1(H;\mathbb Q)=e(H)-v(H)+1$. It follows that 
$$\rank A=e(H)-\dim\ker A\le v(H)-1\le d.$$ By the definition of $A$, $\{s(e):e\in E(H)\}$ lies in an affine subspace of dimension no more than $d-1$. However, by our assumption at least $e(H)\ge q$ distinct elements of $S$ are contained in $\{s(e):e\in E(H)\}$, contradicting the $(d-1,q)$ property of $S$. Therefore $\Gamma_S$ contains no rainbow copy of $H$.

Eventually, note that $\Gamma_S$ is regular with degree $|S|$, we then have
$|S|N^D=e(\Gamma_S)\le\ex^{*}(2N^D,H).$
Taking supremum of $|S|$ over all $(d-1,q)$-general positioned sets $S\subseteq[n]^D$ yields
$$A_{D,d-1,q}(n)\le\frac{\ex^{*}(2N^D,H)}{N^D}.$$
\end{proof}

\end{subsection}

\begin{subsection}{Suspension}\label{sec:graph-suspension}

We next explore how rainbow Tur\'an exponent changes under suspension operations. We will need the following regularizing lemma, where a bipartite graph $F=(A_1\cup A_2, E)$ is said to have \emph{common neighbor property} if for every vertex pair $u_1\ne u_2$ in the same part, there exists $w$ in the other part such that $wu_1,wu_2\in E(F)$. One verifies that, for each bipartite with common neighbor property, the bipartition is unique up to interchanging the two parts.
\begin{lemma}\label{lem:graph-regularize}
Given any properly colored graph $G$, there exists a properly colored bipartite graph $G'$ with the following properties.
\begin{enumerate}[label=\textup{(\roman*)}]
\item $ v(G')\le 8v(G)$, $e(G')\ge \frac 14e(G)$.
\item There exists an integer $\delta$ such that each $v'\in G'$ satisfies $\delta\le d_{G'}(v')\le 2\delta$.
\item For any bipartite $F$ with common neighbor property, if $G'$ has a copy of $F$, then $G$ has an identically colored copy.
\end{enumerate}
\end{lemma}
\begin{proof}
We first find a spanning bipartite subgraph by randomly partitioning $V(G)=V_1\cup V_2$, in the sense that every vertex independently falls into each part with probability $\frac 12$ respectively. Then we obtain a random bipartite subgraph $G_1\subseteq G$ by deleting the edges inside $V_1$ and $V_2$, and we have $\mathbb E(e(G_1))=\frac 12 e(G)$. Therefore there exists a choice of $G_1$ with $e(G_1)\ge \frac 12 e(G)$. Fix one such choice, and let $G_1$ inherit the edge coloring from $G$.

Then we exclude the vertices with small degree. Set $\delta=\ceil{\frac{e(G)}{4v(G)}}$. Iteratively delete vertices $u\in V(G_1)$ from $G_1$ whenever $d_{G_1}(u)<\frac{e(G)}{4(v(G))}$ to obtain a subgraph $G_2\subseteq G_1$. Let $G_2$ inherit the edge coloring from $G_1$. The amount of deleted edges is no more than $v(G_1)\frac{e(G)}{4(v(G))}=\frac 14e(G)$. Therefore $e(G_2)\ge e(G_1)-\frac 14e(G)\ge\frac 14 e(G)$, $v(G_2)\le v(G)$, and $\delta(G_2)\ge \delta$.

Finally, for each vertex of large degree, we create several clones of it and distribute its incident edges among these clones. More precisely, we run Algorithm~\ref{alg:regularize} to split the vertices with large degree and obtain $G'$. 
\begin{algorithm}[!h]
\caption{Splitting large-degree vertices}
\label{alg:regularize}
\renewcommand{\algorithmicrequire}{\textbf{Input:}}
\renewcommand{\algorithmicensure}{\textbf{Output:}}
\begin{algorithmic}[1]
\REQUIRE $G_2$
\ENSURE $G'$
\WHILE{$\exists u\in V(G_2)$ with $d_{G_2}(u)>2\delta$}
\STATE $t\gets\floor{\frac{d_{G_2}(u)}{\delta}}$, $N\gets N_{G_2}(u)$
\STATE Partition $N=N_1\cup\cdots\cup N_t$ with $|N_j|=\delta$, $1\le j<t$ and $\delta\le |N_t|<2\delta$
\STATE $V(G_2)\gets (V(G_2)\backslash\{u\})\cup\{u^{(i)}:1\le i\le t\}$, we call $\{u^{(i)}:1\le i\le t\}$ the clones of $u$.
\STATE $E(G_2)\gets (E(G_2)\backslash\{uw:w\in N\})\cup\{u^{(i)}w:1\le i\le t,w\in N_i\}$
\STATE Let $u^{(i)}w$ inherit the coloring of $uw$.
\ENDWHILE
\RETURN $G_2$
\end{algorithmic}
\end{algorithm}
One easily verifies that the coloring of $G'$ is proper, that $G'$ is also bipartite and all the clones of the same vertex go to the same part, that $e(G_2)$ is preserved at each step, and that the degree of each clone is fixed after the generation step. Therefore $e(G')=e(G_2)\ge \frac14e(G)$, and for each vertex $v'\in V(G')$, $\delta\le d_{G'}(v')\le 2\delta$. It follows that $v(G')\le \frac{2e(G)}{\delta}\le 8v(G)$. 

To verify (\romannumeral3), we show that, if $G'$ contains a copy of $F$, then no two vertices in the copy are clones of the same vertex in $G_2$. Suppose, for the sake of contradiction, that the copy contains clones $u^{(1)}$ and $u^{(2)}$ of $u\in G_2$, then $u^{(1)}$ and $u^{(2)}$ are in the same part of $G'$ and are also in the same part of the copy, since the bipartition of $F$ is unique. By the common neighbor property, there exists $w'$ in the copy with $w'u^{(1)},w'u^{(2)}\in E(G')$. However, the colors of $w'u^{(1)}$ and $w'u^{(2)}$ are inherited from the same edge $uw\in E(G_2)$, where $w'=w$ or $w'$ is a clone of $w$, contradicting the properness of the coloring of $G'$.
\end{proof}

We will also need the following rainbow $C_4$ counting lemma. Write $\#C_4^{rb}(G)$ for the amount of rainbow copies of $C_4$ in an edge-colored graph $G$.

\begin{lemma}\label{lem:graph-rainbow-C4}
For any sufficiently large integer $n$ and any properly colored graph $G$ with $v(G)=n$ and $\delta\le d_G(v)\le2\delta$  for each $v\in V(G)$, if $\delta>100\sqrt{n}$, then $\#C_4^{rb}(G)\ge\frac{1}{100}\delta^4$. 
\end{lemma}
\begin{proof}
Write $\#C_4(G)$ for the amount of $C_4$ copies in $G$, and $\#C_4^{mo}(G)$ for that of $C_4$ copies with at least $2$ edges in the same color.
We estimate $\#C_4(G)$ and $\#C_4^{mo}(G)$ respectively.

For any vertices $v_1\ne v_2$, suppose $u_1\ne u_2$ and $ v_1,v_2\in N_G(u_1)\cap N_G(u_2)$, then $v_1,v_2,u_1$ and $u_2$ form a $C_4$, and each $C_4$ is counted twice in this way. For a fixed pair $v_1,v_2$, the amount of such pairs $u_1,u_2$ is ${{|N_G(v_1)\cap N_G(v_2)|}\choose{2}}$. Therefore, with  $\sum\limits_{v_1,v_2\in V(G)}|N_G(v_1)\cap N_G(v_2)|$ denoted by $S_2$,
$$\#C_4(G)=\frac 12\sum_{v_1,v_2\in V(G)}{{|N_G(v_1)\cap N_G(v_2)|}\choose{2}}\ge\frac{S_2^2}{2n^2-2n}-\frac{S_2}4$$
by Cauchy's inequality.
On the other hand, by double-counting $S_2=\sum\limits_{v\in V(G)}{{d_G(v)}\choose{2}}\in[\frac14n\delta^2,2n\delta^2]$. Therefore $\#C_4(G)\ge\frac{\delta^4}{32}-\frac{\delta^2n}{2}\ge \frac{\delta^4}{40}$ when $n$ is sufficiently large.

Besides, for fixed $v\in V(G)$ and $u_1,u_2\in N_G(v)$, there are at most two choices of $v'$ such that $vu_1v'u_2$ is not rainbow. Indeed, incident edges have different colors by properness, so the possible pairs of edges with the same color are $(vu_1,u_2v')$ and $(vu_2,u_1v')$. Again by properness at $u_2$, there is at most one choice of $v'$ for which $u_2v'$ receives the same color as $vu_1$, and similarly at most one choice of $v'$ makes $(vu_2,u_1v')$ monochromatic. Therefore $\#C_4^{mo}(G)\le2S_2\le4n\delta^2$.
Combining both estimates yields $\#C_4^{rb}(G)\ge\frac{\delta^4}{40}-4n\delta^2\ge\frac{\delta^4}{100}$.
\end{proof}

We are ready to illustrate how rainbow Tur\'an exponent changes under suspension operation.

\begin{proof}[Proof of Theorem \ref{thm:rbturan&suspension}]
Suppose $\ex^{*}(n,F)\le C_0n^{1+\alpha}$. Let $C=2^{26}\cdot100\max\{1,C_0\}$. We prove that, for sufficiently large integer $n$, if graph $G$ satisfies $v(G)=n$ and $e(G)\ge Cn^{1+\frac{1}{2-\alpha}}$, then any proper coloring of $G$ contains a rainbow copy of $\Sigma F$.

We apply Lemma~\ref{lem:graph-regularize} to obtain a regularized graph $G'$. Let $\delta$ be as in Lemma~\ref{lem:graph-regularize}(\romannumeral2), and write $v(G')=n'$. Then $n'\le 8n$, $e(G')\ge\frac14e(G)\ge\frac 14Cn^{1+\frac{1}{2-\alpha}}\ge2^{-8}Cn'^{(1+\frac{1}{2-\alpha})}$, so $\delta\ge\frac{e(G')}{n'}\ge2^{-8}Cn'^{\frac{1}{2-\alpha}}>100\sqrt{n'}$ when $n$ is sufficiently large. Therefore we can apply Lemma \ref{lem:graph-rainbow-C4} to $G'$ to deduce that $G'$ has at least $\frac{\delta^4}{100}$ many rainbow copies of $C_4$.

Now we uniformly and independently assign to each color a label from $\{1,2,3,4\}$, and say a rainbow copy of $C_4$ is \emph{good} if the labels on its four edges, in cyclic order, are $1,2,3,4$ respectively. Each rainbow $C_4$ has probability $4^{-3}$ to be good, so the expected number of good copies of $C_4$ in $G'$ is $4^{-3}\#C_4^{rb}(G')$. Therefore, there is a deterministic way to assign labels, such that the number of good $C_4$ in $G'$ is at least $4^{-3}\#C_4^{rb}(G')$. Furthermore, there exists an edge  $uv$ labeled $1$ admitting at least $\frac{\#C_4^{rb}(G')}{2\cdot4^3n'\delta}$ many good copies of $C_4$, such that every copy contains $uv$, and that in every copy, the edge labeled $2$ is incident to $u$, and that labeled $4$ is incident to $v$.

Therefore we identify a bipartite subgraph $G_{uv}=(N^2(u)\cup N^4(v), E_{uv})$, where $w\in N^2(u)$ if $uw\in E(G')$ and $uw$ is labeled $2$, and similarly, $w\in N^4(v)$ if $vw\in E(G')$ and $vw$ is labeled $4$. For $w_1\in N^2(u)$, $w_2\in N^4(v)$, $w_1w_2\in E_{uv}$ if $w_1w_2\in E(G')$ and $w_1w_2$ is labeled $3$. Note that each edge $w_1w_2\in E_{uv}$ uniquely determines a good $C_4$ containing $uv$, namely $uvw_2w_1$, and vice versa. Therefore $e(G_{uv})\ge \frac{\#C_4^{rb}(G')}{2\cdot4^3n'\delta}\ge\frac{\delta^3}{200\cdot4^3n'}$. Besides, $v(G_{uv})=|N^2(u)|+|N^4(v)|\le d_{G'}(u)+d_{G'}(v)\le 4\delta$. Recall that $\delta\ge 2^{-8}Cn'^{\frac{1}{2-\alpha}}$, we have
$$\frac{e(G_{uv})}{(v(G_{uv}))^{1+\alpha}}\ge\frac{\delta^{2-\alpha}}{100\cdot 4^5n'}\ge\frac{C^{2-\alpha}}{100\cdot2^{26}}>C_0.$$
Hence we can find a rainbow copy of $F$ in $G_{uv}$. This copy of $F$, together with $u$ and $v$, forms a copy of $\Sigma F$ in $G'$. In this copy, for $1\le i\le 4$, denote by $E_i$ the set of edges with label $i$. Then $E_3$ consists of the edges in the copy of $F$, and these edges are in distinct colors. Besides, $E_2$ consists of all the edges at $u$, and the colors of them are distinct by properness at $u$. Similar arguments hold for $E_4$. Since edges with different labels have different colors, we conclude that this copy of $\Sigma F$ is rainbow.

Finally, by Lemma \ref{lem:graph-regularize}(\romannumeral3), note that suspension of any bipartite has the common neighbor property, we can find a copy of $\Sigma F$ in $G$ using the same colors, so this copy is also rainbow, as desired.
\end{proof}

The following corollary indicates how rainbow Tur\'an exponent changes after applying multiple times of suspension operation. As below, the $j$-th suspension of $F$, denoted by $\Sigma^jF$, is iteratively defined by $\Sigma^jF=\Sigma(\Sigma^{j-1}F)$.

\begin{cor}\label{cor:graph-multi-suspension}
Let $F$ be a fixed connected bipartite graph, and suppose $\ex^{*}(n,F)=O(n^{1+\alpha})$ for some $0\le \alpha<1$. For any positive integer $s$, set $\alpha_s=\frac{s-(s-1)\alpha}{s+1-s\alpha}$, then $\ex^{*}(n,\Sigma^sF)=O_s(n^{1+\alpha_s})$.    
\end{cor}

\begin{proof}
Write $\alpha_0=\alpha$, and note that $\alpha_s=\frac{1}{2-\alpha_{s-1}}$ for each positive $s$. We then iteratively apply Theorem \ref{thm:rbturan&suspension} for $s$ times to yield the result.
\end{proof}

\end{subsection}

\begin{subsection}{Quantitative bounds}\label{sec:graph-quant}

To adjust the number of vertices and edges in the optimized graph, we need to add pendent vertices to it, which does not change its rainbow Tur\'an exponent, as shown in the next lemma.

\begin{lemma}\label{lem:graph-pendant}
Let $F$ be a bipartite satisfying $\ex^*(n,F)=O(n^{1+\alpha})$ for some $0\le \alpha<1$, and let $F_{\ell}$ be a graph obtained by arbitrarily adding $\ell$ pendant vertices to $F$, then $\ex^*(n,F_{\ell})=O_{\ell}(n^{1+\alpha}).$    
\end{lemma}

\begin{proof}
Suppose $\ex^*(n,F)\le Cn^{1+\alpha}$ holds for any integer $n$, and $C>e(F)+v(F)+2\ell$. We show that $\ex^*(n,F_{\ell})\le 2Cn^{1+\alpha}$ for any sufficiently large $n$. For any graph $G$ with $v(G)=n$ and $e(G)\ge 2Cn^{1+\alpha}$, iteratively delete vertices $u$ with $d_G(u)<Cn^{\alpha}$ from $G$ to obtain a new graph $G'$. Then $e(G')>e(G)-n\cdot Cn^{\alpha}\ge Cn^{1+\alpha}$, $v(G')\le v(G)$, so $G'$ has a rainbow copy of $F$. To extend this copy by a pendant vertex at some vertex $u$ at each step, we need to avoid at most $e(F)+\ell$ many colors and $v(F)+\ell$ many vertices. This is doable since $d(u)\ge Cn^{\alpha}>e(F)+v(F)+2\ell$, given $n$ is sufficiently large.
\end{proof}

We will start from the rainbow Tur\'an bound for theta graphs $\Theta_{k,t}$ established by Janzer \cite{Janzer2023}. Here $\Theta_{k,t}$ is the union of $t$ internally vertex-disjoint paths of length $k$ with the same endpoints. Hence $\Theta_{k,t}$ is a bipartite graph and we have $v(\Theta_{k,t})=t(k-1)+2$, $e(\Theta_{k,t})=kt$.

\begin{lemma}\label{lem:graph-Theta}
$\operatorname{ex}^{*}(n,\Theta_{k,t})=O_{k,t}(n^{1+\frac{1}{k}}).$
\end{lemma}

We recall the quantitative settings in the definition of $\Gamma(d,q)$. Let $s$ be the unique integer such that $0\le s\le \floor{\frac{d-3}{2}}$ and that $q\in((s+1)(d-s),(s+2)(d-s-1)]$, and put $t=q-(s+1)(d-s)+1$, $k=1+\floor{\frac{d-2s-1}{t}}$. Here we further set $\ell=d-2s-1-t(k-1)$. For fixed $d,q$, we obtain $F$ by arbitrarily adding $\ell$ vertices to $\Theta_{k,t}$, and set $H=\Sigma^s(F)$.

\begin{proof}[Proof of Theorem \ref{thm:graph-main}]
Assume that $q\le\floor{(\frac{d+1}{2})^2}$.
By Theorem \ref{thm:grid&rbturan}, it suffices to show that $v(H)\le d+1$, $e(H)\ge q$, and that $\ex^*(n,H)=O(n^{1+\frac{1}{\Gamma(d,q)}})$. We will show that the first two inequalities are actually equalities. Indeed, $v(F)=v(\Theta_{k,t})+\ell=d-2s+1$, $e(F)=e(\Theta_{k,t})+\ell=d-2s+t-1$, $v(H)=v(F)+2s=d+1$, and $e(H)=e(F)+sv(F)+s^2=(s+1)d-s^2-s+t-1=q$. Finally, combining Lemma \ref{lem:graph-pendant} and Lemma \ref{lem:graph-Theta} yields $\ex^*(n,F)=O(n^{1+\frac 1k})$, and by Corollary~\ref{cor:graph-multi-suspension} applied on $F$ we have $\ex^*(n,H)=O(n^{1+\frac{sk-s+1}{(s+1)k-s}})=O(n^{1+\frac{1}{\Gamma(d,q)}})$, as desired. 
\end{proof}

We remark that, $\Gamma(d,q)$ is the optimal exponent that we can obtain from suspension operations and adding pendant vertices to theta graphs. 

\end{subsection}

\end{section}

\begin{section}{Linear Reduction of tuples}

In Section \ref{sec:fingerprint-vandermonde}, we bound the size of a $(d-1,q)$ general positioned set by assigning each of its tuples of distinct points a shorter representative tuple and show that any representative corresponds to constantly many tuples. Furthermore, in Section \ref{sec:fingerprint-compare}, we compare the graph-based bounds and linear-mapping bounds, and discuss the structural coincidences in specific parameter ranges.

\begin{subsection}{Reduction via full spark linear mappings}\label{sec:fingerprint-vandermonde}
To prove Theorem \ref{thm:fingerprint-main}, we apply a full spark linear mapping on the $h$-tuples of a $(d-1,q)$-general positioned set without repeated coordinates to obtain shorter tuples. We indicate that the amount of possible images is $O(n^{sD})$, and that the preimage at each reduced tuple is bounded to complete the proof. 

More precisely, fix parameters $d,q,h,s,t$ satisfying the quantitative conditions in Theorem \ref{thm:fingerprint-main}, and arbitrarily fix a full spark matrix $\Lambda\in \mathbb Z^{s\times h}$. Then for any integer $1\le j\le h$, any submatrix consisting of $j$ distinct columns has rank $\min\{j,s\}$, and in particular, $\rank(\Lambda)=s$. Moreover, for $S\subseteq\mathbb Z^D$, let $T^h(S)$ be the set of ordered $h$-tuples consisting of distinct elements of $S$.
Then for $\bm x\in T^h(S)$ set $\Phi(\bm x)=\Lambda\bm x\in (\mathbb Z^D)^s$. We also view $\bm x$ as an $(h\times D)$ matrix and $\Phi(\bm x)$ as an $(s\times D)$ matrix.

We first show that, many disjoint tuples consisting of distinct points with the same image would violate the $(d-1,q)$-general position condition.
\begin{lemma}\label{lem:fingerprint-fiber-disjoint}
  Let \(S\subseteq[n]^{D}\) be in \((d-1,q)\) general position with \(\abs{S}\geq q\), then for any full spark linear mapping $\Phi:T^h(S)\rightarrow(\mathbb Z^D)^s$ and any reduced tuple $\bm y\in (\mathbb Z^D)^s$, the preimage of $\bm y$ cannot contain \(t\) tuples, among which no two tuples share any point in $S$.  
\end{lemma}

\begin{proof}
For the sake of contradiction, suppose that $\Phi(\bm x_1)=\Phi(\bm x_2)=\cdots=\Phi(\bm x_t)$, where for each $1\le i\le t$, $\bm x_i=(x_{i,1},\cdots,x_{i,h})$ consists of $h$ distinct points in $S$, and $x_{i_1,j_1}\ne x_{i_2,j_2}$ whenever $i_1\ne i_2$. Set 
$$A=\begin{pmatrix}
-\Lambda & & & & \Lambda\\
& -\Lambda & & & \Lambda\\
& & \ddots & & \vdots\\
& & & -\Lambda & \Lambda
\end{pmatrix}_{s(t-1)\times ht},$$
where each row has $t$ blocks and each column has $t-1$ blocks. The matrix $A$ has full row rank, since the submatrix of the first $t-1$ block columns is block diagonal, and each block is of full row rank. Moreover, set $X=(\bm x_1^{\mathsf T}, \bm x_2^{\mathsf T},\cdots , \bm x_t^{\mathsf T})^{\mathsf T}$, then $AX=0$ since for each $1\le i\le t-1$, $\Lambda \bm x_i=\Lambda \bm x_t$. If we further let $X_1$ be obtained from $X$ by adjoining the all $1$'s vector as the first column, then we also have $AX_1=0$ since the all $1$'s vector clearly lies in the kernel of $A$. By Sylvester's rank inequality, $\rank A+\rank X_1\le ht$, so $\rank X_1\le ht-s(t-1)\le ht-(q-d)$. It follows that $P_0=\{x_{i,j}:1\le i\le t,1\le j\le h\}$, which contains $th$ points, lies in an affine subspace of dimension no more than $th-(q-d)-1$. Arbitrarily add $q-ht$ points in $S\backslash P_0$ to $P_0$ to obtain a new set $P$ of size $q$, then $P$ lies in an affine subspace of dimension no more than $d-1$, contradicting the $(d-1,q)$ general position condition.
\end{proof}

We then proceed to uniformly bound the size of the preimages at any reduced tuple of $\Phi$.

\begin{lemma}\label{lem:fingerprint-fiber-size}
There exists a constant $C=C(d,q,s,h,t)$ such that for any \(S\subseteq[n]^{D}\) in \((d-1,q)\) general position with \(\abs{S}\geq q\), any full spark linear mapping $\Phi:T^h(S)\rightarrow(\mathbb Z^D)^s$, and any reduced tuple $Y\in (\mathbb Z^D)^s$, the preimage of $Y$ contains at most $C$ distinct elements of $T^h(S)$.    
\end{lemma}

\begin{proof}
For each fixed tuple $(q,d,s,t)$, we prove by induction on $h$ for $h\in[s,\frac qt]$. When $h=s$, the full spark matrix $\Lambda$ is invertible, and $\Phi$ is injective. Therefore the preimage of each reduced tuple contains at most one element.    

If the conclusion holds for $h=k$, where $s\le k\le\frac qt-1$, we consider the case $h=k+1$ as follows. Fix $Y=(\bm y_1,\cdots,\bm y_s)^{\mathsf T}\in (\mathbb Z^D)^s$ and denote by $F_Y$ the preimage of $Y$. Let $\bm x_1,\cdots,\bm x_{t_0}$ be a maximal disjoint family of tuples in $F_Y$, by Lemma \ref{lem:fingerprint-fiber-disjoint} we have $t_0<t$. Let $S_Y=\{x_1,x_2,\cdots,x_{(k+1)t_0}\}$ be all the points in $S$ that appear in some $\bm x_i$ ($1\le i\le t_0$). By maximality of $\bm x_1,\cdots,\bm x_{t_0}$, every tuple $\bm x\in F_Y$ contains some point in $S_Y$. We then have $F_Y=\bigcup\limits_{i=1}^{(k+1)t_0}\bigcup\limits_{j=1}^{k+1}F_{i,j}$, where $F_{i,j}$ consists of all the tuples in $F_Y$ whose $j$-th coordinate is $x_i$. 

For fixed $1\le i\le (k+1)t_0$, $1\le j\le k+1$, delete the $j$-th column of $\Lambda$ to obtain a submatrix $\Lambda_j$, and obtain $F_{i,j}'$ from $F_{i,j}$ by deleting the $j$-th coordinate from each tuple. Furthermore, set $\Phi_j(\bm x')=\Lambda_j\bm x'$ for any $\bm x'\in F'_{i,j}$, then $\Phi_j$ is also full spark, and $\Phi_j(\bm x')=(\bm y_1-\Lambda_{1j}x_i,\bm y_2-\Lambda_{2j}x_i,\cdots,\bm y_s-\Lambda_{sj}x_i)^{\mathsf T}$. That is, all tuples of $F'_{i,j}$ have the same images under $\Phi_j$. Let $C_0$ be the constant from inductive hypothesis, then $|F_{i,j}|=|F'_{i,j}|\le C_0$. This holds for every $1\le i\le(k+1)t_0$
and $1\le j\le k+1$. Hence
\[|F_Y|\le\sum_{i=1}^{(k+1)t_0}\sum_{j=1}^{k+1}|F_{i,j}|\le(k+1)^2t_0C_0<(k+1)^2tC_0.\]
This holds for each $Y\in (\mathbb Z^D)^s$, so we complete the proof of the case $h=k+1$, and by induction we are done.
\end{proof}

Furthermore, we upper bound the size of the image of $\Phi$ and deduce a bound of $|T^h(S)|$ for each $S\subseteq[n]^D$ in $(d-1,q)$-general position, from which we yield an upper bound of $|S|$.

\begin{proof}[Proof of Theorem \ref{thm:fingerprint-main}]
For each $S\subseteq[n]^D$ in $(d-1,q)$  general position, take $\Lambda\in \mathbb Z^{s\times h}$ to be any full spark matrix, and denote the maximum absolute value of its entries by $K$. We also assume that $|S|\ge q$. Since each entry of $\Lambda$ has absolute value no larger than $K$, we have $|\Phi(T^h(S))|\le (2hKn)^{sD}$. Let $C$ be the constant from Lemma \ref{lem:fingerprint-fiber-size}, then $|T^h(S)|\le C|\Phi(T^h(S))|=O(n^{sD})$. On the other hand, $|T^h(S)|=\prod\limits_{i=0}^{h-1}(|S|-i)=\Omega(|S|^h)$, thus we conclude that $|S|=O(n^{\frac{sD}{h}})$. This holds for every $S$ in $(d-1,q)$ general position, so $A_{D,d-1,q}(n)=O(n^{\frac{sD}{h}})$, as desired.  

Finally, when $s<h$, we use Theorem \ref{thm:HD&grid} to show that $\HD_d(p,q)\ge p^{\frac hs-o(1)}$. If $h=s$, we use the trivial bound $\HD_d(p,q)\ge p-q$.
\end{proof}

We then prove Corollary~\ref{cor:fingerprint-exp} by properly choosing the parameters.

\begin{proof}[Proof of Corollary \ref{cor:fingerprint-exp}]
For each $s\in\{1,2,\cdots,d-1\}$, pick $t=1+\ceil{\frac{q-d}{s}}$, $h=s+\floor{\frac{d-s}{t}}\ge s$, $q'=s(t-1)+d\ge q$, then $th=ts+t\floor{\frac{d-s}{t}}\leq st+d-s=q'$, and $s(t-1)=q'-d$. Therefore by Theorem \ref{thm:fingerprint-main} we have $\HD(p,q')\ge p^{\frac hs-o(1)}$. By monotonicity we have $\HD(p,q)\ge\HD(p,q')\ge p^{\frac hs-o(1)}$. Taking the maximum over $s\in\{1,2,\cdots,d-1\}$ gives us the desired result. 
\end{proof}
\end{subsection}    

\begin{subsection}{Comparison of graph-based and linear reduction methodology}\label{sec:fingerprint-compare}
By comparing the exact value of $\Gamma(d,q)$ and $\Pi(d,q)$, we note the following coincidences.   

\begin{prop}\label{prop:fingerprint-identical}
For any $d\ge 2$,
\begin{enumerate}[(a)]
\item If $q=(r+1)(d-r)$ for some $1\le r\le \floor{\frac{d-1}{2}}$, then $\Gamma(d,q)=\Pi(d,q)=\frac{r+1}{r},$
\item If $d<q<d+\sqrt d$, then $\Gamma(d,q)=\Pi(d,q)=\floor{\frac{q}{q-d+1}},$
\item If $q>\floor{\frac{(d+1)^2}{4}}$, then $\Gamma(d,q)=\Pi(d,q)=1$.
\end{enumerate}
\end{prop}

In fact, these quantitative coincidences arise from the underlying structural agreement, as the following analysis reveals.

Recall that, in both derivations we show that some forbidden structures could provide too many independent affine relations, contradicting the $(d-1,q)$-general position condition. More precisely, in the graph-based estimation, the affine relations come from the cycle space of a bipartite whose rainbow copies are forbidden, while in the linear reduction method, these relations are provided by the row space of 
$$A=\begin{pmatrix}
-\Lambda & & & & \Lambda\\
& -\Lambda & & & \Lambda\\
& & \ddots & & \vdots\\
& & & -\Lambda & \Lambda
\end{pmatrix}_{s(t-1)\times ht},$$
where $\Lambda$ is the full spark reduction matrix. We will show that, in the cases (a)(b) of Proposition \ref{prop:fingerprint-identical}, by properly choosing $\Lambda$, these two spaces are identical. To avoid confusion of notation, we write $s_G$ and $t_G$ for the parameters $s$ and $t$ in the graph-based construction, denoting the number of suspensions and the number of paths in the theta graph respectively, while writing $s_F$ and $t_F$ for the parameters $s$ and $t$ in the linear mapping construction.

In the cases of (a), for the graph construction we pick $k=2$, $t_G=d-2r+1$, $s_G=r-1$, $\ell=0$, so $\Theta_{k,t_G}=K_{2,d-2r+1}$, and $H=\Sigma^{r-1}\Theta_{k,t_G}=K_{r+1,d-r}$. For the linear reduction method, the parameters are $s_F=r$, $h=r+1$, $t_F=d-r$, and we particularly pick 
$$\Lambda=\begin{pmatrix}
-1 & & & & 1\\
& -1 & & & 1\\
& & \ddots & & \vdots\\
& & & -1 & 1
\end{pmatrix}_{r\times (r+1)}.$$
We now identify the column indices of \(A\) with the edges of \(H\). Let the two parts of $H$ be $V=\{v_1,\cdots,v_{r+1}\}$ and $U=\{u_1,\cdots,u_{d-r}\}$, and identify the $i$-th column of the $j$-th block column of $A$ with the edge $v_iu_j$ for each $i\in[r+1]$, $j\in[d-r]$. Then the $i$-th row vector of the $j$-th block row has $-1$ at $v_iu_j$, $1$ at $v_{r+1}u_j$, $1$ at $v_iu_{d-r}$, $-1$ at $v_{r+1}u_{d-r}$, and $0$ at other coordinates, which is precisely the cycle vector of the $4$-cycle $v_iu_jv_{r+1}u_{d-r}$, up to a global swap of sign. Moreover, by fixing a spanning forest $F\subset H$, where $E(F)=\{v_iu_{d-r}:i\in[r+1]\}\cup\{v_{r+1}u_j:j\in[d-r]\}$, the 4-cycle $v_iu_jv_{r+1}u_{d-r}$ is the unique cycle in $F\cup\{v_iu_j\}$. By Lemma \ref{lem:graph-cycle-space}, all these cycle vectors form a basis of the cycle space of $H$. Therefore the row space of $A$ and the cycle space of $H$ are identical.

In the cases of (b), for the graph construction we pick $k=\floor{\frac{q}{q-d+1}}$, $t_G=q-d+1$, $\ell=q-kt$ and $s_G=0$, so the graph $H$ is $\Theta_{k,t_G}$ attached with $\ell$ pendant vertices. For the linear reduction method, the parameters are $s_F=1$, $h=\floor{\frac{q}{q-d+1}}=k$, $t_F=q-d+1$. This time we let $\Lambda$ be the row vector $(1,-1,\cdots,(-1)^{h-1})_{1\times h}$, then $A$ has $ht_F=kt_G$ many columns. Denote the $i$-th edge of the $j$-th path of $\Theta_{k,t_G}$ by $e_{i,j}$, and identify it with the $i$-th column of the $j$-th block column of $A$, where $i\in[k]$, $j\in[t_G]$. Then for each $i\in[t_G-1]$, the $i$-th row vector has $1$ and $-1$ alternating on the $i$-th and the $t_G$-th path, which is precisely the cycle vector of the $2k$-cycle consisting of these two paths. Again these cycle vectors make up the cycle space of $H$, thus the row space of $A$ and the cycle space of $H$ are identical.

Furthermore, Proposition \ref{prop:fingerprint-identical}(c) indicates a common limitation of the two methods, which is the consequence of the following elementary result.

\begin{lemma}\label{lem:fingerprint-largest-bip}
Any bipartite $F$ with $v(F)-c(F)\le d$ satisfies $e(F)\le\floor{(\frac{d+1}{2})^2}$.   
\end{lemma}

\begin{proof}
Identifying one vertex from each connected component makes $F$ connected without changing the value $v(F)-c(F)$ and $e(F)$. So we may assume that $F$ is connected and that $v(F)\le d+1$ at the beginning. Suppose the bipartition is $V(F)=V_1\cup V_2$ with $|V_1|=d'$, then $e(F)\le d'(d+1-d')\le\floor{(\frac{d+1}{2})^2}$.
\end{proof}

Therefore, in the graph-based construction, for any fixed $d$, the optimized bipartite to which we apply Theorem \ref{thm:grid&rbturan} has no more than $\floor{(\frac{d+1}{2})^2}$ edges, so by the quantitative restriction, the method works only when $q\le\floor{(\frac{d+1}{2})^2} $. 

In the linear reduction method, we also use a bipartite-structured construction to show that, once $q>\floor{(\frac{d+1}{2})^2}$, the size of the preimages of any non-trivial reduction mapping cannot be uniformly bounded. More precisely, given any integer $N$, for $D>N+q$ and $h>s$, any $s\times h$ matrix $\Lambda$ admits a $(d-1,q)$ general positioned set $S\subseteq [n]^D$ such that $\Phi:T^h(S)\rightarrow (\mathbb Z^D)^s$, $\Phi(\bm x)=\Lambda \bm x$ with at least $N$ distinct tuples mapped to the same image by $\Phi$. Indeed, pick $a_i=e_i$ for each $1\le i\le h$, $b_i=e_{h+i}$ for each $1\le i\le N$, where $e_i$ is the $i$-th basis vector of $\mathbb Z^D$ with $1$ on the $i$-th coordinate and $0$ on other coordinates, pick a non-zero vector $\bm w=(w_1,\cdots,w_h)\in \ker \Lambda$, and set $S=\{a_i+w_ib_j:i\in [h],j\in[N]\}$. Up to a translation we assume that $S\subseteq[n]^D$ for sufficiently large $n$. Furthermore, set $\bm x_j=(a_1+w_1b_j,\cdots,a_h+w_hb_j)\in T^h(S)$ for each $j\in[N]$, then $\Phi(\bm x_1)=\cdots=\Phi(\bm x_N)$ since $\Lambda w=0$. It remains to show that $S$ is in $(d-1,q)$ general position.

\begin{lemma}\label{lem:fingerprint-explode-fiber}
Any $(d-1)$-dimensional affine subspace intersects $S$ at not more than $\floor{(\frac{d+1}{2})^2}$ points.   
\end{lemma}

\begin{proof}
Without loss of generality, suppose for some $1<h_0\le h$, $w_i\ne 0$ if and only if $i\in[h_0]$.
For any $T\subseteq S$, construct an auxiliary bipartite $G_T'$ on $[h_0]\cup [N]$ by adding an edge between $i$ and $j$ whenever $a_i+w_ib_j\in T$, and obtain $G_T$ from $G_T'$ by adding an isolated edge for each $a_i\in T$, $i>h_0$.
We show that $\dim \aff(T)=v(G_T)-c(G_T)-1$.

Partition $T=T_1\cup T_2$, where $T_1=\{a_i+w_ib_j\in T:i\in [h_0]\}$ corresponds to $G_T'$, $T_2=\{a_i\in T:h_0<i\le h\}$ corresponds to the matching part of $G_T$. By definition each element in $T_2$ is linearly independent from other elements in $T$, so $\dim \aff(T)=\dim \aff(T_1)+|T_2|$.
Besides, each element of $T_2$ contributes $2$ to $v(G_T)$ and $1$ to $c(G_T)$, so $v(G_T)-c(G_T)-1=v(G_{T}')-c(G_T')-1+|T_2|$. Hence it suffices to show that $\dim \aff(T_1)=v(G_T')-c(G_T')-1$.

Let $Y$ be the matrix whose column vectors are precisely the elements of $T_1$, then the column for $(a_i+w_ib_j)$ has $1$ at the $i$-th row, $w_i$ at the $(h+j)$-th row and $0$ at other positions. By multiplying the column for $(a_i+w_ib_j)$ by $\frac{1}{w_i}$ and the $i$-th row by $w_i$ for each $i\in[h_0]$, the submatrix consisting of the first $h_0$ rows and the $h+1$-th to the $h+N$ rows of $Y$ aligns with the incidence matrix of $G_T'$, say $M$, and all other rows are empty. By Lemma \ref{lem:graph-cycle-space}(a), $\rank Y=\rank M=v(G_T')-c(G_T')$. It follows that $\dim \mathrm{span}(T_1)=\rank Y=v(G_T')-c(G_T')$. Finally, all the elements of $Y$ lie in the affine subspace $\sum\limits_{i=1}^{h_0}x_i=1$, which does not contain the origin, thus $\dim \aff(T_1)=\dim \mathrm{span}(T_1)-1=v(G_T')-c(G_T')-1$.

Therefore $\dim \aff(T)=v(G_T)-c(G_T)-1$ for each $T\subseteq S$. If $\dim \aff(T)\le d-1$, then $v(G_T)-c(G_T)\le d$, and by Lemma \ref{lem:fingerprint-largest-bip} and the construction of $G_T$, $|T|=e(G_T)\le \floor{(\frac{d+1}{2})^2}$, as desired.
\end{proof}

Since $q>\floor{(\frac{d+1}{2})^2}$, any $q$ points of $S$ do not lie in an affine subspace of dimension $d-1$, as desired. We conclude that, the limitations of both methods are somehow the consequences of Lemma \ref{lem:fingerprint-largest-bip}.

\end{subsection}

\end{section}

\begin{section}{Concluding Remark}
In this paper, we establish improved lower bounds for Hadwiger–Debrunner numbers using hypergraph containers, rainbow Turán theory and linear reduction of tuples. To complement these results, we conclude with the following remarks.

We first remark that, our three constructions are complementary, and for each method, there are ranges for $d,q$ in which it yields a strictly greater exponent for $p$ than both alternatives, illustrated as follows.

\begin{prop}\label{prop:remark-quant}
Write $\Psi(d,q)=\frac{d+q-2}{q-1}$, which is the exponent of $p$ in the bound provided by hyperplane construction. Then
\begin{enumerate}[(a)]
\item When $q>\floor{(\frac{d+1}{2})^2}$, $\Psi(d,q)=\frac{d+q-2}{q-1}>1=\Gamma(d,q)=\Pi(d,q)$, so hyperplane construction provides strictly better bounds.
\item For integer $a\ge 3$, when $d=2a+1$, $q=4a+1$, we have $\Gamma(d,q)>\max\{\Psi(d,q),\Pi(d,q)\}$, so the graph-based construction provides strictly better bounds.
\item For integer $a\ge 4$, when $d=3a+2$, $q=5a$, we have $\Pi(d,q)>\max\{\Psi(d,q),\Gamma(d,q)\}$, so the linear mapping construction provides strictly better bounds.
\end{enumerate}
\end{prop}

Again, we write $s_G$ for the parameter $s$ in the formulation of $\Gamma(d,q)$.

\begin{proof}[Proof of Proposition \ref{prop:remark-quant}]
We only prove (b) and (c). When $d=2a+1$, $q=4a+1$ for $a\ge 3$, we have $\Psi(d,q)=\frac 32$. In the formulation of $\Gamma(d,q)$, we have $s_G=1$, $t=2$, $k=a$, so $\Gamma(d,q)=2-\frac 1a>\Psi(d,q)$. To see $\Gamma(d,q)>\Pi(d,q)$, by definition it suffices to show that, for each $1\le s\le d-1$, $\frac 1s\floor{\frac{2a+1-s}{1+\ceil{2a/s}}}<1-\frac{1}{a}$. When $s=1$, $\floor{\frac{2a+1-s}{1+\ceil{2a/s}}}=0$, so the inequality clearly holds. When $s\ge 2$, $\frac 1s\floor{\frac{2a+1-s}{1+\ceil{2a/s}}}\le\frac{2a+1-s}{2a+s}\le1-\frac{3}{2a+2}<1-\frac 1a$. We complete the proof of (b).

When $d=3a+2$, $q=5a$ for $a\ge 4$, we have $\Psi(d,q)=1+\frac{3a+1}{5a-1}<2$. For the graph construction, we have $s_G=0$, $t=2a-1$, $k=2$, so $\Gamma(d,q)=2$. For the linear reduction construction, by taking $s=2$ we have $\Pi(d,q)\ge 1+\frac 12\floor{\frac{d-2}{1+\ceil{(q-d)/2}}}=\frac 52>\max\{\Psi(d,q),\Gamma(d,q)\}$. We complete the proof of (c).
\end{proof}

We finally remark on the common limitation of the three methods. Notably, in all the three methods, we use the maximum depth, namely the maximum amount of members that share a common point, to bound the piercing number. The following result in \cite{KellerSmorodinskyTardos2018} yields a limitation of this strategy.
\begin{theorem}[\cite{KellerSmorodinskyTardos2018}, Proposition 2.3]\label{thm:remark-depth}
Let $\cF$ be a family of $n\ge 2p$ convex sets in $\mathbb R^d$ with the $(p,q)$-property. Then there exists a point that pierces $\Omega_d(p^{-\frac{q-1}{q-d}}qn)$ elements of $\cF$.
\end{theorem}

Therefore, in any finite compact convex family with the $(p,q)$-property, the ratio between its size and its maximum depth is $O_{d,q}(p^{\frac{q-1}{q-d}})$, and the lower bound of Hadwiger-Debrunner exponent provided by the maximum depth argument is no greater than $\frac{q-1}{q-d}$. To achieve a further breakthrough, one has to discover the families with the $(p,q)$-property and non-uniform intersection patterns simultaneously, and employ better strategies to control the piercing numbers.
\end{section}

\subsection*{Acknowledgement}
The authors would like to thank Zixiang Xu for useful suggestions on writing and a series of valuable discussions on combinatorial convexity. Wenchong Chen would like to thank Cosmin Pohoata for insightful discussions on relative topics.

\subsection*{Declaration of the use of AI}
GPT-6 astra was used to assist in developing the proof of Theorem~\ref{thm:fingerprint-main}, identifying a suitable hypergraph container theorem, as well as optimizing the quantitative results and preparing the
manuscript. All AI-suggested arguments were independently verified, revised and streamlined by
the authors, who take full responsibility for the correctness of all mathematical results presented in
this paper.

\bibliographystyle{plain}
\bibliography{bib}

\appendix
\end{document}